\input psfig.sty

\overfullrule=0pt

{\centerline {\bf SMOOTH COMBS INSIDE HEDGEHOGS} }

\bigskip
\bigskip

{\centerline {KINGSHOOK BISWAS}}

\medskip

{\centerline {Department of Mathematics}}

{\centerline {University of California, Los Angeles}}

{\centerline {Los Angeles, CA 90095}}

\bigskip
\bigskip

{\bf Abstract.} {\it We use techniques of tube-log Riemann surfaces due to R.Perez-Marco to construct a hedgehog
containing smooth $C^{\infty}$ combs. The hedgehog is a common hedgehog for a family of commuting non-linearisable
holomorphic maps with a common indifferent fixed point. The comb is made up of smooth curves, and is transversally
bi-H\"older regular. }

\bigskip


\medskip


\bigskip

\noindent {\bf 1. Introduction.}

\bigskip

We consider the dynamics of a holomorphic diffeomorphism $f(z) = e^{2\pi i \alpha} z + O(z^2), \, \alpha \in
\bf{R} - \bf{Q}$, defined in a neighbourhood of the indifferent irrational fixed point 0. The map $f$ is said to
be {\it linearisable} if there is a holomorphic change of variables $z = h(w) = w + O(w^2)$ such that $$ h^{-1}
\circ f \circ h = R_\alpha $$ in a neighbourhood of the origin, where $R_\alpha(w) = e^{2\pi i \alpha} w$ is the
rigid rotation.

\medskip

The problem of {\it linearisation}, or determining when such an $f$ is linearisable, has a long and interesting
history, including the work of H. Cremer (see [2], [3]) in the 1920's, C.L.Siegel (see [10]) in 1942, A.D.Bruno in
the 1960's, and Yoccoz (see [11]), who resolved the problem of the optimal arithmetic condition on $\alpha$ for
linearisability in 1987. In this article however we are primarily concerned with the structure of invariant sets
for the dynamics near the fixed point, rather than the issue of linearisability.

\medskip

When $f$ is linearisable, $f$ behaves like the rotation by angle $2\pi\alpha$ around 0 on the maximal domain of
linearisation called the {\it Siegel disk} of $f$, which turns out to be a domain of the form $h( \{ |w| < R \})$
where $R$ is less than or equal to the radius of convergence of $h$ around 0. In this case the compacts $h( \{ |w|
\leq r \})$  for $r < R$ are clearly completely invariant under $f$. However, even in the more general case where
$f$ is not necessarily linearisable, R.Perez-Marco found completely invariant continua for $f$ which persist:

\medskip

\noindent {\bf Theorem (Perez-Marco, [4]).} {\it Let $f(z) = e^{2\pi i \alpha}z + O(z^2), \, \alpha \in \bf{R} -
\bf{Q}$ be such that $f$ and $f^{-1}$ are defined and univalent on a neighbourhood of the closure of a Jordan
domain $\Omega \subset \bf{C}$ containing $0$. Then there exists a full compact connected set $K$ contained in
$\overline{\Omega}$ such that $0 \in K$, $K \cap \partial\Omega \neq \phi$ and $f(K) = f^{-1}(K) = K$. }

\medskip

The invariant compacts $K$ thus obtained are called {\it Siegel compacta}. A classical topological theorem of
G.D.Birkhoff (see [1]) does in fact guarantee, for planar homeomorphisms near Lyapunov unstable points, the
existence of compacts, but which are either positive or negative invariant, not necessarily totally invariant. In
the holomorphic setting one can thus improve to obtain totally invariant compacts. Moreover, though this will not
be used in what follows, if the boundary of the Jordan domain $\Omega$ is $C^1$ smooth, then Perez-Marco has shown
that such a compact $K$ is in fact unique (see [8]). If $K$ is not contained in the closure of a linearisation
domain it is called a {\it hedgehog}. A {\it hedgehog} is called linearisable or non-linearisable depending on
whether it contains a linearisation domain or not. We will be mostly concerned with non-linearizable hedgehogs in
this paper. Perez-Marco has shown that these have no interior (see [8]). The structure of such hedgehogs is
topologically complex; for example, Perez-Marco shows in [5] that they are not locally connected at any point
different from the fixed point. Objects such as {\it combs} (homeomorphs of the product of a Cantor set and an
interval) had been expected to be found within hedgehogs. Indeed this is stated by Perez-Marco in [7], where he
also sketches how the techniques of "tube-log Riemann surfaces" used in [7] may be adapted to construct a hedgehog
containing combs. We carry out this suggested construction in this article. The comb obtained is quite regular,
being made up of smooth curves containing a dense set of analytic curves. Thus hedgehogs can exhibit some
smoothness. We have the following:

\medskip

\noindent {\bf Main Theorem.} {\it There exists a Cantor set $C \subset \bf{R}$ and a family of commuting
holomorphic maps $(f_t(z) = e^{2\pi i t}z + O(z^2))_{t\in C}$, all defined and univalent on a common neighbourhood
of the closed unit disk $\overline{\bf{D}}$, which have for irrational times $t$ a common hedgehog $K \subset
\overline{\bf{D}}$, such that $K$ contains a comb $\cal{E}$, the homeomorphic image of $C_0 \times [1,2]$ for an
explicit Cantor set $C_0$. The comb is made up of smooth $C^{\infty}$ curves, and is bi-H\"older regular. More
precisely, for the homeomorphism

$$\eqalign{ \Psi : & C_0 \times [1,2] \to \cal{E} \cr
       & (\theta \ \, , \ \ t) \quad \mapsto \Psi(\theta,t) \cr
}$$
we have:

\medskip

(a) For each fixed $\theta$, the map $t \mapsto \Psi(\theta,t)$ is smooth.

\medskip

(b) $\Psi$ is $\alpha$-H\"older regular for all $\alpha, \, 0 < \alpha < 1$, and ${\Psi}^{-1}$ is Lipschitz.

\medskip

The maps $f_t$ are non-linearisable for $t \in C \cap (\bf{R} - \bf{Q})$. }

\medskip

We note that since the hedgehog is non-linearisable, the comb
obtained is non-trivial, in the sense that it is not contained
inside the interior of the hedgehog. We first describe the
construction, explaining how the hedgehog and comb are formed, and
then carry out the proofs of the properties listed above.

\medskip

For other examples of hedgehogs, we refer the reader to the articles of Perez-Marco [6, 7], in the first of which
non-linearisable hedgehogs are constructed along with an explicit coding of the dynamics on the hedgehog, while in
the second linearisable hedgehogs containing Siegel disks with smooth boundary are constructed.

\bigskip

\noindent {\bf 2.  The construction.}

\bigskip

\noindent {\bf 2.1   General idea of the construction.}

\medskip

The construction we will present is based on the techniques of "tube-log Riemann surfaces" invented by Perez-Marco
in [9], where they were first used. They were used again in [6] and [7]. It is strongly advised to read this
section in conjunction with [6] and [7]. As these techniques are explained in detail in [6] and especially in [7],
and the basic elements of the construction are the same as in [7], to avoid repetition the description below is
shorter and slightly informal.

\medskip

We start with the linear flow of translations $(F_t(z)=z+t)_{t\in \bf{R}}$  on the upper half-plane $\bf{H}$ ({\bf
note} : for convenience, instead of working with dynamics and dynamical objects such as Siegel disks, hedgehogs
etc defined on the unit disk ${\bf D}$, we will actually work instead with their lifts to the upper half-plane via
the universal covering $E : {\bf H} \to {\bf D} - \{0\}, E(z) = e^{2\pi i z})$.

\medskip

\noindent {\bf 2.1.1)   The Riemann surface ${\cal S}_0$.}

\medskip

The idea is 'fold' the dynamics by introducing nonlinearities in the dynamics, such as fixed points; these result
in a shrinking of the linearisation domain. This is accomplished by means of the following Riemann surface with
distinguished charts ${\cal S}_0$ (see Perez-Marco [7] section 1.2.1) :

\medskip

{\hfill {\centerline {\psfig {figure=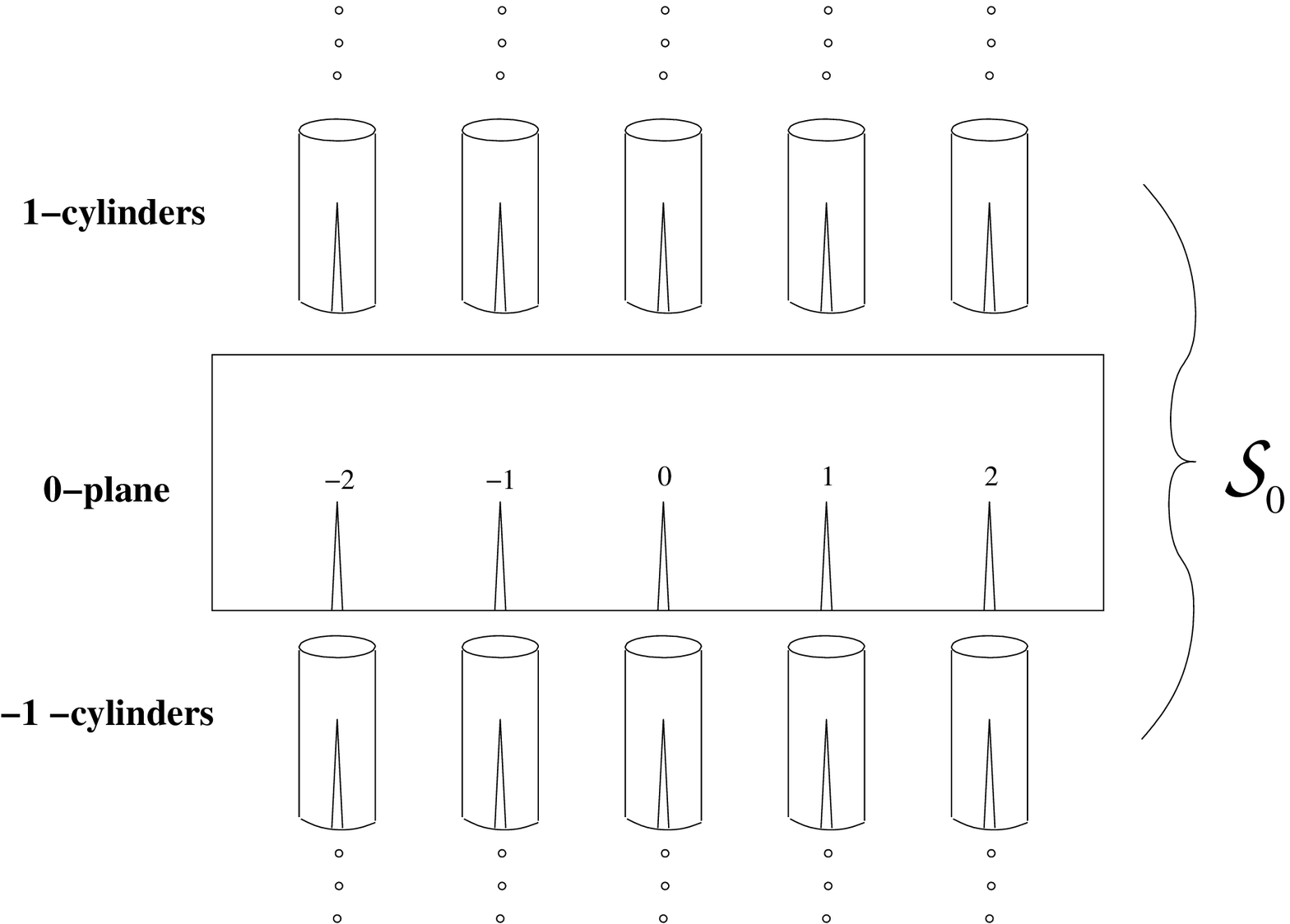,height=8cm}}}}

{\centerline {\bf Figure 1}}

\medskip

The surface ${\cal S}_0$ is formed from a copy of $\bf{C}$ with 'slits' $[n, n-i\infty]$ at each of the integers
$n\in \bf{Z}$, above and below each of which is pasted isometrically along the slits a family (indexed by
$\bf{Z}$) of slit cylinders ${\bf{C}}/\bf{Z}$ with slits $[0,0-i\infty]$. This construction gives not only a
Riemann surface, but also a canonical set of charts for the surface, which allow us to write formulas for
functions defined on the surface in terms of these canonical coordinates. We note here that the surface ${\cal
S}_0$ has 'ramification points of infinite order' at the points corresponding to $\bf Z$ in the $0$-plane and to
$0$ in the cylinders. These can be seen as points added in the completion of the flat metric on ${\cal S}_0$
(induced from the flat metrics on its building blocks). However it is important to note that these points do not
belong to the surface ${\cal S}_0$ itself, nor can the complex structure on ${\cal S}_0$ be extended to these
points.

\medskip

\noindent {\bf 2.1.2)   Uniformization of the Riemann surface ${\cal S}_0$.}

\medskip

The upper ends of the cylinders on the other hand do correspond to points for the complex structure. Adding in
these points gives a simply connected Riemann surface which it is then not too hard to see is biholomorphic to
$\bf{C}$. So ${\cal S}_0$ should be biholomorphic to the complex plane minus a doubly infinite discrete set of
points; in fact, there is an explicit formula for the uniformization (see Perez-Marco [7], section 1.2.3), which
in terms of the distinguished charts on ${\cal S}_0$ is given by $$ K(z) = {1 \over 2\pi i} \log ( {-\log
({1-e^{2\pi iz}})}) $$

The map $K$ satisfies the properties $$\displaylines{ K(z+1) = K(z)+1 \quad \hbox{ on all of the 0-sheet of ${\cal
S}_0$ } \cr K(z)-z \to 0 \quad \hbox{ as Im } z \to + \infty \hbox{ in the 0-sheet of ${\cal S}_0$} \cr }$$ for
the appropriate choices of the branches of the logarithms.

\medskip

The images under $K$ of horizontal lines just above the real line $\bf{R}$ in the 0-plane are thus 1-periodic
(with respect to parametrization by the $x$-coordinate) curves in $\bf{C}$; since the ramification points at $\bf
Z$ in the $0$-plane correspond to ends at infinity under the uniformization, these curves oscillate with
amplitudes which tend to infinity as we take lines tending to the real line. The upper ends of the cylinders
correspond under $K$ to a discrete set of points $P \subset \bf{C}$; amongst these, those that correspond to the
upper ends of the 1-cylinders lie below, but close to, the image of the real line. The figure below illustrates
the image of the foliation of horizontal lines on ${\cal S}_0$ under $K$.

\medskip

{\hfill {\centerline {\psfig {figure=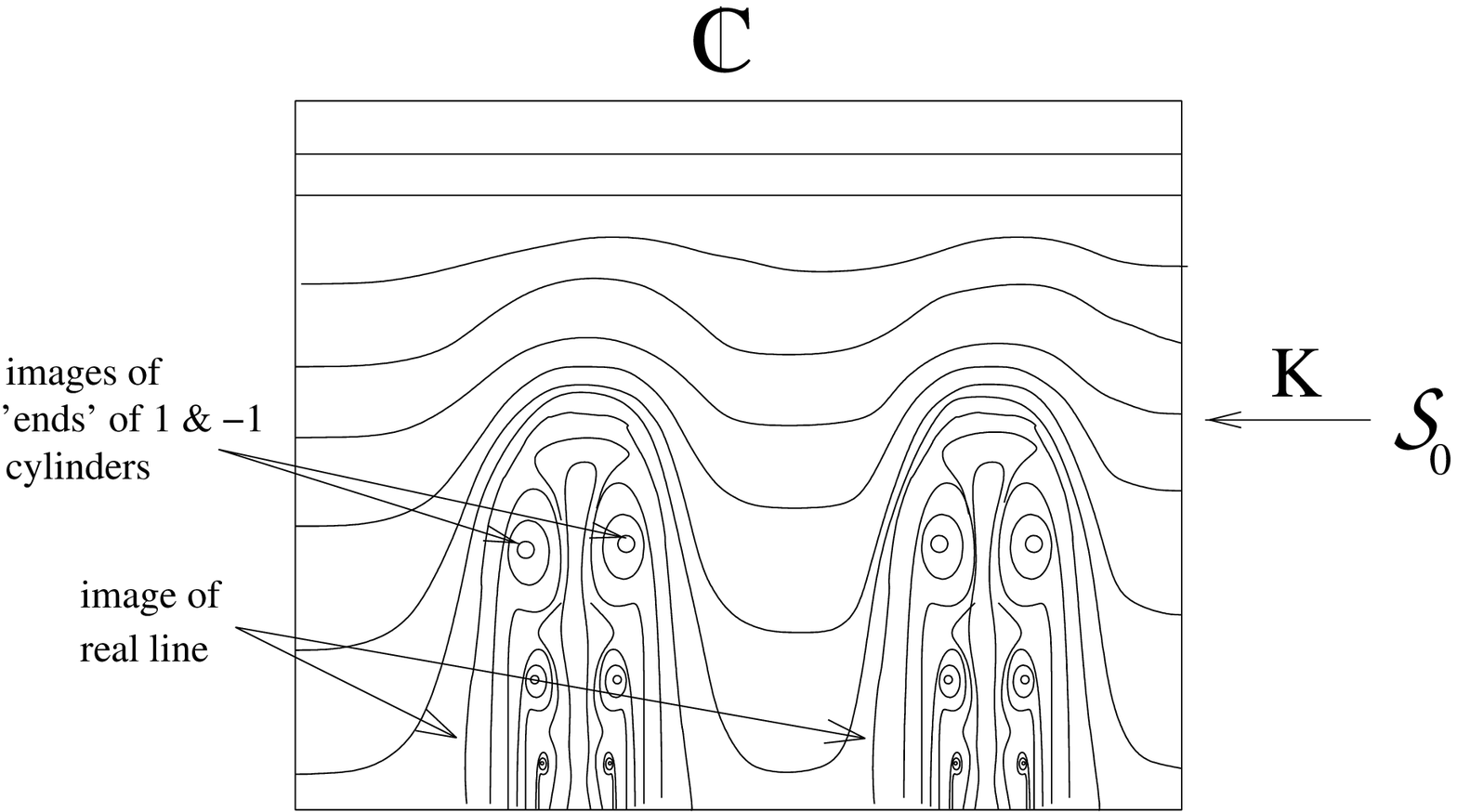,height=6cm}}}}

{\centerline {\bf Figure 2}}

\medskip

\noindent {\bf 2.1.3)   Lifting and folding of dynamics.}

\medskip

Consider the linear flow of translations $(F_t(z)=z+t)_{t\in \bf{R}}$  on the upper half-plane $\bf{H}$. The upper
half-plane of the $0$-plane of ${\cal S}_0$ projects univalently onto the upper half-plane $\bf{H}$, so
 we can lift each $F_t$ to a map $\tilde{F}_t(z)$ from the upper half-plane of the $0$-plane of ${\cal S}_0$ to
itself, given in terms of the distinguished chart by the same formula $\tilde{F}_t(z) = z+t$. We would like to
extend the lift analytically  to a univalent map $\tilde{F}_t$ of the surface to itself; however, it is not
possible to extend to the points ${\bf Z} - t$ of the $0$-plane or the points $0 - t$ of the cylinders, as these
would have to be mapped to the ramification points, and moreover there is a non-trivial monodromy when continuing
analytically around these points.

\medskip

To overcome this problem, we can remove from the surface a small $\epsilon$-neighbourhood (for the flat metric on
${\cal S}_0$) ${\cal I}(\epsilon) \subset {\cal S}_0$ of the ramification points. On ${\cal S}_0 - {\cal
I}(\epsilon)$, for $|t| < \epsilon$, the aforementioned problem does not arise for the maps $\tilde{F}_t$, each of
which then extends to a univalent map $\tilde{F}_t : {\cal S}_0 - {\cal I}(\epsilon) \to {\cal S}_0$, which is
expressed in terms of the distinguished charts on ${\cal S}_0$ by the same formula $\tilde{F}_t(z) = z+t$ (it
suffices in fact to remove only a set of horizontal slits ${\bf Z} + (\epsilon, \epsilon)$ near the ramification
points, but for consistency with the rest of the article we consider removing neighbourhoods).

\medskip

Let $T : {\cal S}_0 \to {\cal S}_0$ be the automorphism given in the charts on ${\cal S}_0$ by $T(z) = z+1$. We
note that the lifts of the integer translations $z \mapsto z + n$ extend to all of ${\cal S}_0$ as the
automorphisms $T^n : {\cal S}_0 \to {\cal S}_0$. These automorphisms commute with the maps $\tilde{F}_t$, hence we
can also define for times $t' = n + t \in {\bf Z} + (-\epsilon, \epsilon)$ the lifts $\tilde{F}_{t'} : {\cal S}_0
- {\cal I}(\epsilon) \to {\cal S}_0$ by $\tilde{F}_{n + t} := T^n \circ \tilde{F}_t$. Hence we obtain a semi-flow
$(\tilde{F}_t)$ defined for times $t \in {\bf Z} + (-\epsilon,\epsilon)$ on a large part ${\cal S}_0 - {\cal
I}(\epsilon)$ of the surface.

\medskip

The 'upper ends' of the cylinders on the surface are points for the complex structure, which become fixed points
(modulo the automorphism $T : {\cal S}_0 \to {\cal S}_0$) for the dynamics on the surface. Hence lifting the
dynamics to the surface has the effect not only of extending its domain of definition, but also of introducing
nonlinearities, such as fixed points.

\medskip

To see this dynamics on the surface as dynamics on the plane, we use the uniformization $K : {\cal S}_0 \to {\bf
C} - P$ :

\medskip

Conjugating the lifts $\tilde{F}_t$ by $K$ we obtain dynamics $K \circ \tilde{F}_t \circ K^{-1}$ on a large part
$K({\cal S}_0 - {\cal I}(\epsilon)) \subset {\bf C}$ of the plane. In particular any half-plane $\{ \hbox{ Im  } z
\geq -y \} \subset {\bf C}$ where $y \geq 0$ is contained in $K({\cal S}_0 - {\cal I}(\epsilon))$ for $\epsilon$
small enough (see [7], section 1.2.4), so we can always ensure a large domain of definition for the dynamics. The
dynamics flows along the foliation shown in Figure 2, and has as fixed points (modulo $\bf Z$) the points $P$
corresponding under $K$ to the upper ends of the cylinders. This dynamics is what we mean by the 'folded'
dynamics.

\medskip

\noindent {\bf 2.1.4)   Normalization of the uniformization.}

\medskip

The images under $K$ of horizontal lines in the upper half-plane of the $0$-sheet of ${\cal S}_0$ are graphs over
the $x$-axis which are 1-periodic with respect to parametrization by the $x$-coordinate. For the image of a
horizontal line at a height $\delta > 0$ in the $0$-sheet, its peaks and troughs occur at the points $K(n+1/2 +
i\delta)$ and $K(n + i\delta)$ respectively (where $n$ runs over all integers). The amplitude of oscillation (the
difference in height between peaks and troughs) is a strictly decreasing function of $\delta$ which tends to
$+\infty$ as $\delta$ tends to $0$, and tends to $0$ as $\delta$ tends to $+\infty$.

\medskip

We define a family of normalized uniformizations as follows:

\medskip

First we place the ramification points of ${\cal S}_0$ slightly below the real axis instead of on the real axis;
more precisely, for $\delta > 0$ we denote by ${\cal S}_{\delta}$ the Riemann surface obtained by the same
construction as ${\cal S}_0$, pasting a plane and cylinders isometrically along slits, but with ramification
points placed at ${\bf Z} - i\delta$ in the $0$-plane and $0 - i\delta$ in the cylinders. The uniformization of
this surface is given by the map  $z \mapsto K(z+i\delta)$ instead of $z \mapsto K(z)$. The image under this
uniformization of the real line ${\bf R}$ in the $0$-sheet of ${\cal S}_{\delta}$ is a 1-periodic graph with an
amplitude of oscillation that is large when $\delta$ is small. For small $\delta$, the troughs of this curve have
large negative imaginary parts and lie well below the real axis; in applications however, we would like the image
of the real line to lie above the $x$-axis, with the troughs touching the $x$-axis, and so we define the
uniformization $K_{\delta} : {\cal S}_{\delta} \to {\bf C}$ of ${\cal S}_{\delta}$ by $$ K_{\delta}(z) :=
K(z+i\delta) - i \times \hbox{"height of the troughs"} $$ where by the "height of the troughs" we mean $$
\hbox{"height of the troughs"} \ = \hbox{ Im } K(n + i\delta) \ = - {1 \over 2\pi } \log (-\log (1-e^{-2\pi
\delta}))$$ (where $n$ here is any integer). The curve $K_{\delta}({\bf R})$ then lies in the upper half-plane of
${\bf C}$ with troughs lying on the real axis at the points $K_{\delta}(n + 0i ) = n + 0i, \, n \in {\bf Z}$; the
image of the upper half-plane of the $0$-sheet of ${\cal S}_{\delta}$ under $K_{\delta}$, which we denote by
$K_{\delta}({\bf H})$ is the region in the upper half-plane of ${\bf C}$ bounded by $K_{\delta}({\bf R})$.

\medskip

In applications it will also be useful to have the image of the real line oscillate at higher frequencies, more
precisely with period $1/a$ instead of period $1$ for integers $a \geq 1$; to do this we consider the
uniformizations ${1 \over a} K_{\delta}$ obtained by dividing $K_{\delta}$ by integers $a \geq 1$. The image of
the real line,  the curve ${1 \over a} K_{\delta}({\bf R})$, is then a $1/a$-periodic curve. Finally, in order to
control the amplitude of oscillation, we note that given $a \geq 1$ integer and $h > 2$ real, there exists a
unique $\delta = \delta(a,h)$ (as is shown in Perez-Marco [7] section 1.2.5), with $0 < \delta < 1/2$, such that
the amplitude of oscillation of $K_{\delta}({\bf R})$ is $a \cdot h$, and hence that of ${1 \over a}
K_{\delta}({\bf R})$ is $h$. Moreover $\delta(a,h) \to 0$ for fixed $h$ when $a \to \infty$. The constant $\delta$
satisfies (see Perez-Marco [7] section 1.2.5) $$ {\log \left(1 + e^{-2\pi \delta} \right) \over -\log \left(1 -
e^{-2\pi \delta} \right)} = e^{-2\pi ah} $$ The peaks and troughs of the curve ${1 \over a} K_{\delta}({\bf R})$
occur at the points ${1 \over a} K_{\delta}(n+1/2 + 0i) = n+1/2 + ih$ and ${1 \over a} K_{\delta}(n + 0i) = n +
0i$ respectively, and we have $\mathop{\rm Max}_{x \in {\bf R}} $ Im ${1 \over a} K_{\delta}(x) = $ Im ${1 \over
a} K_{\delta}(n+1/2) = h$, \ $\mathop{\rm Min}_{x \in {\bf R}} $ Im ${1 \over a} K_{\delta}(x) = $ Im ${1 \over a}
K_{\delta}(n) = 0$, for any integer $n \in {\bf Z}$.

\medskip

{\hfill {\centerline {\psfig {figure=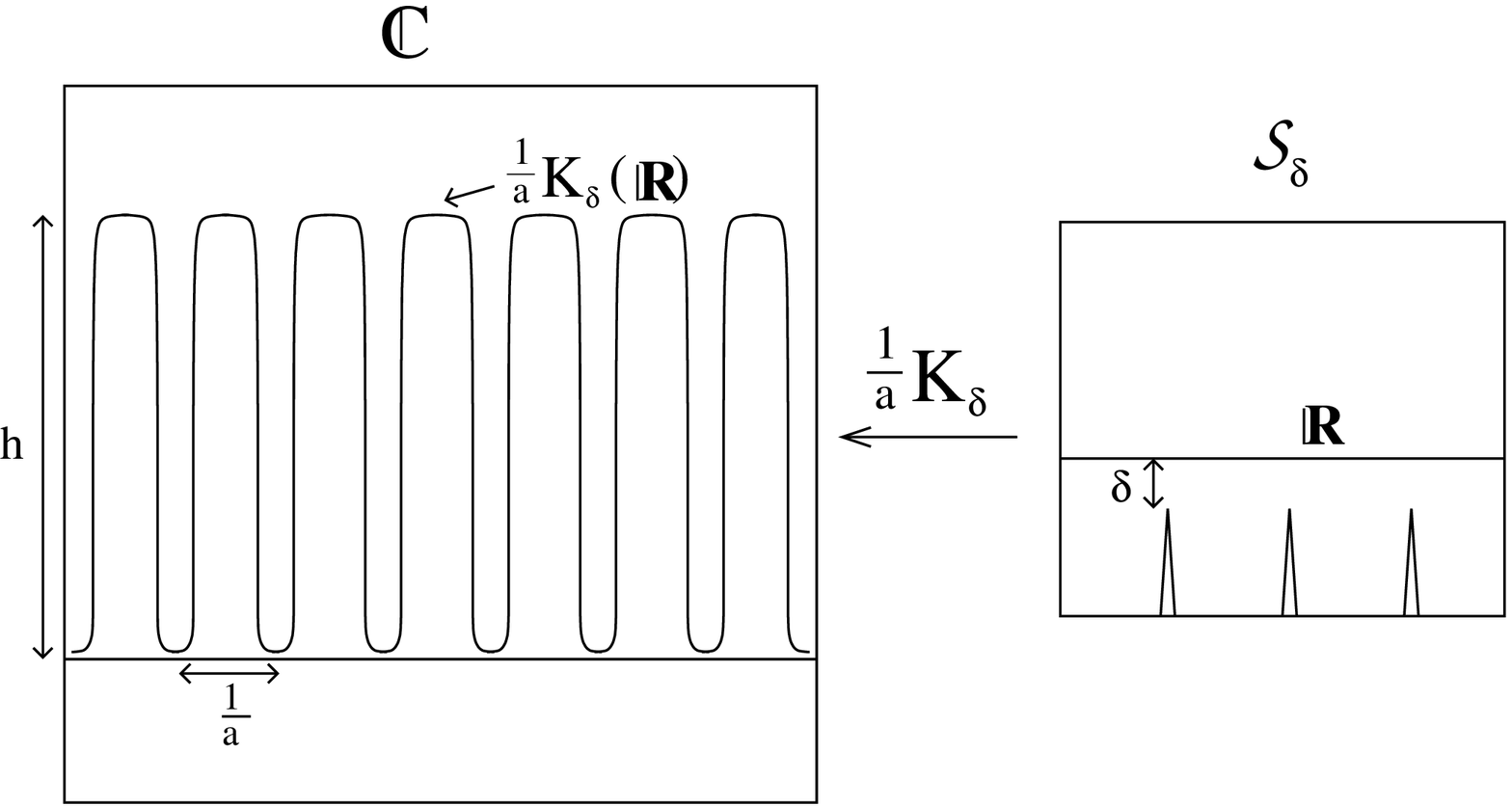,height=7cm}}}}

{\centerline {\bf Figure 3.}}

\medskip

We thus have a family of normalized uniformizations ${1 \over a} K_{\delta}$ depending on two parameters $a \geq
1$ an integer and $h > 2$ a real, where here $\delta$ is understood to be $\delta = \delta(a,h)$. The maps ${1
\over a} K_{\delta}$ have the following properties :

\medskip

\noindent {\bf (1)} ${1 \over a} K_{\delta}$ is univalent on a neighbourhood of $\overline{ \bf{H} }$, mapping
$\bf{H}$ into itself.

\medskip

\noindent {\bf (2)} ${1 \over a} K_{\delta}(z + 1) = {1 \over a} K_{\delta}(z) + {1 \over a}$

\medskip

\noindent {\bf (3)} ${1 \over a} K_{\delta}(i {\bf R}_{\geq 0}) = i {\bf R}_{\geq 0}, \, ({1 \over a}
K_{\delta})'(z) > 0$ for $z \in i {\bf R}_{\geq 0}$

\medskip

\noindent {\bf (4)} There exists a constant $C$ independent of $a$ and $h > 2$ such that $$ \left|{1 \over a}
K_{\delta}(z) - \left({1 \over a}z + ih \right) \right| \leq {C \over a} $$ for all $z$ with Im $z > 1$.

\medskip

\noindent {\bf (5)} Let $z_0 \in {\bf C}$ denote the limit of ${1 \over a} K_{\delta}(z)$ when $z \to \infty$ in
${\cal S}_{\delta}$ through the upper end of one of the 1-cylinders. There exists a constant $C_1$ independent of
$a$ and $h > 2$ such that $$ \hbox{Im } z_0 \geq h - {C_1 \over a} $$ (note that the images under ${1 \over a}
K_{\delta}$ of the upper ends of the 1-cylinders differ from each other by $1/a$-translations, so Im $z_0$ is
independent of the choice of 1-cylinder). We note that, for the dynamics obtained by conjugation by ${1 \over a}
K_{\delta}$, these points are no longer necessarily fixed points (modulo $\bf{Z}$), but rather periodic points
(modulo $\bf{Z}$) of period $a$.

\medskip

We note that the properties {\bf (1)-(4)} refer to ${1 \over a} K_{\delta}(z)$ for $z$ in the $0$-plane of ${\cal
S}_{\delta}$, while {\bf (5)} refers to the 1-cylinders. These properties appear in Perez-Marco [7], section
1.2.5, but can also be checked by elementary computations . We note that condition {\bf (5)} allows us to control
the heights in the upper half-plane of the non-linearities of the dynamics obtained by conjugation by ${1 \over a}
K_{\delta}$; indeed {\bf (5)} says that the heights of the periodic orbits is approximately $h$, hence larger
choices of $h$ result in somehow smaller linearisation domains.

\medskip

\noindent {\bf 2.1.5)   Construction of the limit dynamics.}

\medskip

The main idea is, roughly, to iterate the process described in section 2.1.3 above. We take a surface ${\cal
S}_{\delta_1}$, lift the linear flow of translations to ${\cal S}_{\delta_1}$ and conjugate by a uniformization
$K_0 = {1 \over a_0} K_{\delta_1}$ to obtain folded dynamics as described in section 2.1.3. Having fixed ${\cal
S}_{\delta_1}$ and $K_0$, we take a surface ${\cal S}_{\delta_2}$, lift the linear flow of translations to ${\cal
S}_{\delta_2}$ and conjugate by a uniformization $K_1 = {1 \over a_1} K_{\delta_2}$; we then lift this folded
dynamics to ${\cal S}_{\delta_1}$ and conjugate by the uniformization $K_0$, to obtain a twice-folded dynamics,
given by conjugating by the map $K_0 K_1$. We continue this way, choosing at each stage a surface ${\cal
S}_{\delta_n}$ and a uniformization $K_{n-1}$, to obtain at stage $n$ a dynamics given by conjugation by $K_0
\dots K_{n-1}$. At each stage the dynamics obtained is linearisable, but the linearisation domain decreases; the
point is that, by choosing the appropriate uniformizations and controlling the geometry of the domains obtained at
each stage, we should be able to obtain, in the limit, dynamics having an invariant domain with the desired
geometry.

\medskip

In the construction that follows, we define a sequence $$K_n = {1 \over a_n} K_{\delta_{n+1}}, \ \delta_{n+1} =
\delta(a_n, h_n), \ n \geq 0 $$ of uniformizations depending on two sequences of parameters $(a_n)_{n \geq 0}$ and
$(h_n)_{n \geq 0}$. The whole construction is completely determined by the sequences of parameters $a_n$ and
$h_n$, which we will indicate how to choose in the course of the construction. For the moment we assume that they
are given, and construct the desired limit dynamics as follows:

\medskip

Let $n \geq 1$.

\medskip

We let ${\cal F}_{n,n} = (F_{n,n,t}(z) = z + t)_{t \in {\bf R}}$ be the linear flow of translations. As explained
in section 2.1.3 above, we can lift the elements of this flow to univalent maps $\tilde{F}_{n,n,t}$ defined on the
part of the surface ${\cal S}_{\delta_{n}}$ obtained by deleting an $\epsilon$-neighbourhood ${\cal I}(\epsilon)$
of the ramification points and for times $t \in {\bf Z} + (-\epsilon, \epsilon)$. We take $\epsilon =
3{a_n}^{-1}$, and conjugate by the uniformization $K_{n-1}$ of ${\cal S}_{\delta_{n}}$ to obtain a semi-flow
${\cal F}_{n-1,n} = (F_{n-1,n,t})_{t \in A_{n-1,n}}$ defined by $$ F_{n-1,n,{a_{n-1}}^{-1}t} := K_{n-1} \circ
\tilde{F}_{n,n,t} \circ K_{n-1}^{-1} $$ where $$ {a_{n-1}}^{-1}t \in A_{n-1,n} := {a_{n-1}}^{-1} ({\bf Z} +
[-2{a_n}^{-1}, 2{a_n}^{-1}]). $$ (we could in this step have taken a bigger set of times for $A_{n-1,n}$ such as
${a_{n-1}}^{-1} ({\bf Z} + (-3{a_n}^{-1}, 3{a_n}^{-1}))$, but in the steps that follow it will be necessary to
restrict to smaller sets of times). We note that the time has been reparametrized in order to have the right
behaviour at $+i\infty$, when Im $z \to +\infty$, $$ F_{n-1,n,t}(z) = z + t + o(1). $$ These maps are defined on
$K_{n-1}({\cal S}_{\delta_{n}} - {\cal I}(3{a_n}^{-1})) \subset {\bf C}$, which contains a large half-plane ${\bf
H}_{-y} = \{ \hbox{ Im } z \geq -y \}$, where $y \to +\infty$ as $a_n \to +\infty$.

\medskip

We would now like to lift the semi-flow ${\cal F}_{n-1,n}$ to the surface ${\cal S}_{\delta_{n-1}}$ and conjugate
by the uniformization $K_{n-2}$ to obtain a semi-flow ${\cal F}_{n-2,n}$. The difference with the previous case is
that we are now starting with a nonlinear flow instead of the linear flow of translations. However, by giving up
some space in the lower region of the half-plane of definition ${\bf H}_{-y}$, we can work with maps close to
linear ones; more precisely we have (see [7] Lemma 1.1.6 and Proposition 2.2.6), for $t \in A_{n-1,n}$, if Im $z
\geq -y + {1 \over 2\pi} \log a_{n-1} + C$, (where $C$ is a universal constant) then $$ |F_{n-1,n,t}(z) - (z + t)|
< {1 \over a_{n-1}} $$ Since $y \to +\infty$ as $a_n \to +\infty$, we observe that once $a_{n-1}$ is fixed, by
choosing $a_n$ large enough we can ensure that this inequality holds in as large a half-plane ${\bf
H}_{-\tilde{y}}$ as desired, where $-\tilde{y} = -y + {1 \over 2\pi} \log a_{n-1} + C$.

\medskip

Now consider the region ${\cal J}(-\tilde{y}) = \{ z \in {\cal S}_{\delta_{n-1}} : \hbox{ Im }z \leq -\tilde{y} \}
\subset {\cal S}_{\delta_{n-1}}$ (note that there is a well-defined function Im $z$ on ${\cal S}_{\delta_{n-1}}$).
Since the estimate above holds on the complement ${\cal S}_{\delta_{n-1}} - {\cal J}(-\tilde{y})$ of this region,
it should be clear that if we now remove as well an $\epsilon$-neighbourhood ${\cal I}(\epsilon)$ of the
ramification points, this time with $\epsilon = 3{a_{n-1}}^{-1}$, then we can define, for $t \in {\bf Z} +
A_{n-1,n} \cap [-2{a_{n-1}}^{-1}, 2{a_{n-1}}^{-1}]$, the lifts of elements of ${\cal F}_{n-1,n}$ as univalent maps
$\tilde{F}_{n-1,n,t} : {\cal S}_{\delta_{n-1}} - ({\cal J}(-\tilde{y}) \cup {\cal I}(3{a_{n-1}}^{-1})) \to {\cal
S}_{\delta_{n-1}}$.

\medskip

We can now proceed as before and define the semi-flow ${\cal F}_{n-2,n} = (F_{n-2,n,t})_{t \in A_{n-2,n}}$ by $$
F_{n-2,n,{a_{n-2}}^{-1}t} := K_{n-2} \circ \tilde{F}_{n-1,n,t} \circ K_{n-2}^{-1} $$ where $$ {a_{n-2}}^{-1}t \in
A_{n-2,n} := {a_{n-2}}^{-1} ({\bf Z} + A_{n-1,n} \cap [-2{a_{n-1}}^{-1}, 2{a_{n-1}}^{-1}]) $$ (the times have
again been reparametrized appropriately so that $F_{n-2,n,t}(z) = z+t+o(1)$ when Im $z \to +\infty$).

\medskip

These maps are defined on $K_{n-2}\left({\cal S}_{\delta_{n-1}} - ({\cal J}(-\tilde{y}) \cup {\cal
I}(3{a_{n-1}}^{-1}))\right) \subset {\bf C}$, which contains a large half-plane ${\bf H}_{-y'} = \{ \hbox{ Im } z
\geq -y' \}$, where $y' \to +\infty$ as $a_{n-1} \to +\infty$ and $\tilde{y} \to \infty$ (see [7] Proposition
2.2.8).

\medskip

Following this same procedure of lifting and conjugating the dynamics, we can define successively the semi-flows
${\cal F}_{n-3,n}, \, {\cal F}_{n-4,n}, \dots, \, {\cal F}_{0,n}$, with ${\cal F}_{j,n} = (F_{j,n,t})_{t \in
A_{j,n}}$ defined for times $t \in A_{j,n}$ where $A_{j,n}$ is defined inductively by $$ A_{j,n} := {a_{j}}^{-1}
({\bf Z} + A_{j+1,n} \cap [-2{a_{j+1}}^{-1}, 2{a_{j+1}}^{-1}]) $$ for $j=n-1,\dots,0$ and $A_{n,n} = {\bf R}$. We
note here that the sets $A_{0,n}$ form a decreasing sequence, ie $A_{0,n+1} \subset A_{0,n}, n \geq 1$; indeed, $$
A_{n-1,n+1} = {a_{n-1}}^{-1} ({\bf Z} + A_{n,n+1} \cap [-2{a_n}^{-1}, 2{a_n}^{-1}]) \subset {a_{n-1}}^{-1} ({\bf
Z} + [-2{a_n}^{-1}, 2{a_n}^{-1}]) = A_{n-1,n} $$ from which it follows similarly that $$\eqalign{ A_{n-2,n+1} &
\subset A_{n-2,n} \cr A_{n-3, n+1} & \subset A_{n-3,n} \cr & \dots \cr A_{0,n+1} & \subset A_{0,n} \cr }$$.

Thus we obtain a sequence of semi-flows ${\cal F}_{0,n} = (F_{0,n,t})_{t \in A_{0,n}}$, where each one is given by
'folding' the linear flow $n$ times. Moreover, the $a_n$'s can be chosen inductively to grow fast enough so that
all the semi-flows ${\cal F}_{0,n}$ are defined in the upper half-plane ${\bf H}$, or indeed in any fixed
half-plane ${\bf H}_{-M} = \{ \hbox{ Im } z > -M \}$, as large as desired; indeed we can state as a proposition
the following (which follows directly from Perez-Marco [7] sections 1.3.1 and 2.2.3):

\medskip

\noindent {\bf Proposition 2.1.1.} \ {\it Given $M > 0$, there exists a sequence of conditions $({\cal C}_n)_{n
\geq 1}$ on the integers $a_n$ of the form $$ ({\cal C}_n) \qquad \quad a_n \geq C_n(a_0, \dots, a_{n-1}, h_0,
\dots, h_{n-1}, M) $$ such that when they are fulfilled then all the semi-flows $({\cal F}_{0,n})_{n \geq 1}$ are
defined in the half-plane ${\bf H}_{-M} = \{ \hbox{ Im } z > -M \}$.}

\medskip

We observe thus a feature of the construction which is crucial when passing to the limit, namely that it allows us
to successively paste nonlinearities to the dynamics without decreasing its domain of definition too much,
allowing us to retain a half-plane ${\bf H}_{-M}$.

\medskip

Assume that the $a_n$'s have been chosen so that all the semi-flows ${\cal F}_{0,n}$ are defined on ${\bf
H}_{-M}$. From the sequence of semi-flows ${\cal F}_{0,n}$ we can extract a convergent subsequence (it always
exists, see Perez-Marco [7], section 1.1.2) converging normally for times $t \in A = \cap_{n \geq 0} A_{0,n}$ to
get a limit semi-flow ${\cal F} = (F_t)_{t \in A}$ defined on the half-plane ${\bf H}_{-M}$ for times $t \in A$.
We note that elements of the limit semi-flow do indeed commute, since we can pass to the limit in the equation
$$F_{0,n_k,t} \circ F_{0,n_k,s} (z) = F_{0,n_k,s} \circ F_{0,n_k,t} (z) $$ (using uniform convergence in compact
neighbourhoods of the points $z, F_t(z), F_s(z)$).

\medskip

The Cantor set of times $C$ in the main theorem is given by $$\displaystyle{ C := \bigcap_{n \geq 0} A_{0,n} /
{\bf Z} }$$ If all the $a_n$'s are large enough, say for example $a_n \geq 5$ for all $n$, then the set $C$ is
indeed a Cantor set.

\medskip

The semi-flow $(f_t)_{t \in C}$ of the main theorem is defined on a common neighbourhood $E({\bf H}_{-M})$ of the
closed unit disk $\overline{{\bf D}}$ via the equation $$ f_t \circ E = E \circ F_t, \ t \in C $$ (where $E$ is
the universal covering $E : {\bf C} \to {\bf C} - \{0\}, E(z) = e^{2\pi i z}$).

\medskip

As stated in section 2.1.4, each uniformization $K_n$ maps the upper half-plane of the $0$-plane of ${\cal
S}_{\delta_{n+1}}$ into the upper half-plane ${\bf H}$. So we can consider its restriction to the upper-half plane
of the $0$-plane as a map $K_n : {\bf H} \to {\bf H}$ of ${\bf H}$ into itself. If we consider the composition of
these maps $K_0 \circ K_1 \circ \dots \circ K_n$, then it is clear that the set $K_0 K_1 \dots K_n (\overline{{\bf
H}})$ is totally invariant under the semi-flow ${\cal F}_{0,n+1}$, and hence the set $${\cal H} = \bigcap_{n =
0}^{\infty} K_0 K_1 \dots K_n (\overline{{\bf H}}) $$ is then totally invariant under the limit semi-flow ${\cal
F}$. Moreover, ${\cal H}$, being the decreasing limit of the domains $K_0 K_1 \dots K_n (\overline{{\bf H}})$, is
full, and ${\cal H} \cap \partial {\bf H} \neq \phi$ (since $0 \in {\cal H}$). Thus ${\cal H}$ is in fact a common
hedgehog for irrational elements of the limit semi-flow ${\cal F}$. The hedgehog $K$ of the semi-flow $(f_t)_{t
\in C}$ in the main theorem is given by $$ K = E({\cal H}) $$

\medskip

It is important to note in this construction that the only requirement imposed on the $a_n$'s so that the
limit dynamics is defined on the upper half-plane is that they grow fast enough, as expressed by Proposition
2.1.1 above. Subject to this restriction, the construction still offers enough flexibility to construct
hedgehogs with quite different geometries by varying the choices of the parameters $a_n$ and $h_n$.

\medskip

In the following sections we show how to choose the uniformizations $K_n$ appropriately (by proper choices of the
parameters $a_n, h_n$), and control the geometries of the domains $K_0 K_1 \dots K_n(\overline{{\bf H}})$, so that
the hedgehog ${\cal H}$ contains a comb with the desired geometry.

\bigskip

\noindent {\bf 2.2  Formation of the comb.}

\medskip

The purpose of this section is to describe informally the idea behind
the construction of the comb within the hedgehog; precise
statements and proofs are given in the following section, 2.3.

\medskip

We define $$\displaylines{ \hbox{The strip } S := \{ 1 \leq \hbox{ Im } z \leq 2 \} \cr \hbox{The vertical
half-lines } \xi_x := \{ x+iy : y \geq 0 \} \, , \, x \in {\bf R} \cr }$$ In this section we prove only the
following Proposition:

\medskip

\noindent {\bf Proposition 2.2.1} \ {\it Let $m \in {\bf Z}$ and $N \geq 0$ be integers. Then the curve $K_0 \dots
K_N(\xi_m)$ is contained in the hedgehog ${\cal H}$.}

\medskip

\noindent {\bf Proof :} The properties {\bf (2)} and {\bf (3)} of the maps ${1 \over a} K_{\delta}$ given in
section 2.1.4 imply that for any integer $m \in {\bf Z}$  and any $n \geq 0$ we have $$ K_n(\xi_m) = \xi_{m/a_n}.
$$ Thus if we consider a curve $K_0 \dots K_N(\xi_m)$, then for all $n > N$ we can write $$\eqalign{ K_0 \dots
K_N(\xi_m) & = K_0 \dots K_N \dots K_n \left(\xi_{m \cdot a_{N+1} \dots a_n } \right) \cr
                     & \subset K_0 \dots K_N \dots K_n (\overline{{\bf H}}) \cr
}$$ and hence $$ K_0 \dots K_N(\xi_m) \subset \bigcap_{n = N+1}^{\infty} K_0 \dots K_n(\overline{{\bf H}}) = {\cal
H} $$ (since the domains in the intersection are decreasing, the intersection taken from $N+1$ to $\infty$ is the
same as from $1$ to $\infty$). $\diamondsuit$

\medskip

Thus all curves of the form $K_0 \dots K_N(\xi_m)$ are contained in the hedgehog ${\cal H}$. It is with such
curves that we will construct the comb. We consider the formation of the comb in stages:

\medskip

{\hfill {\centerline {\psfig {figure=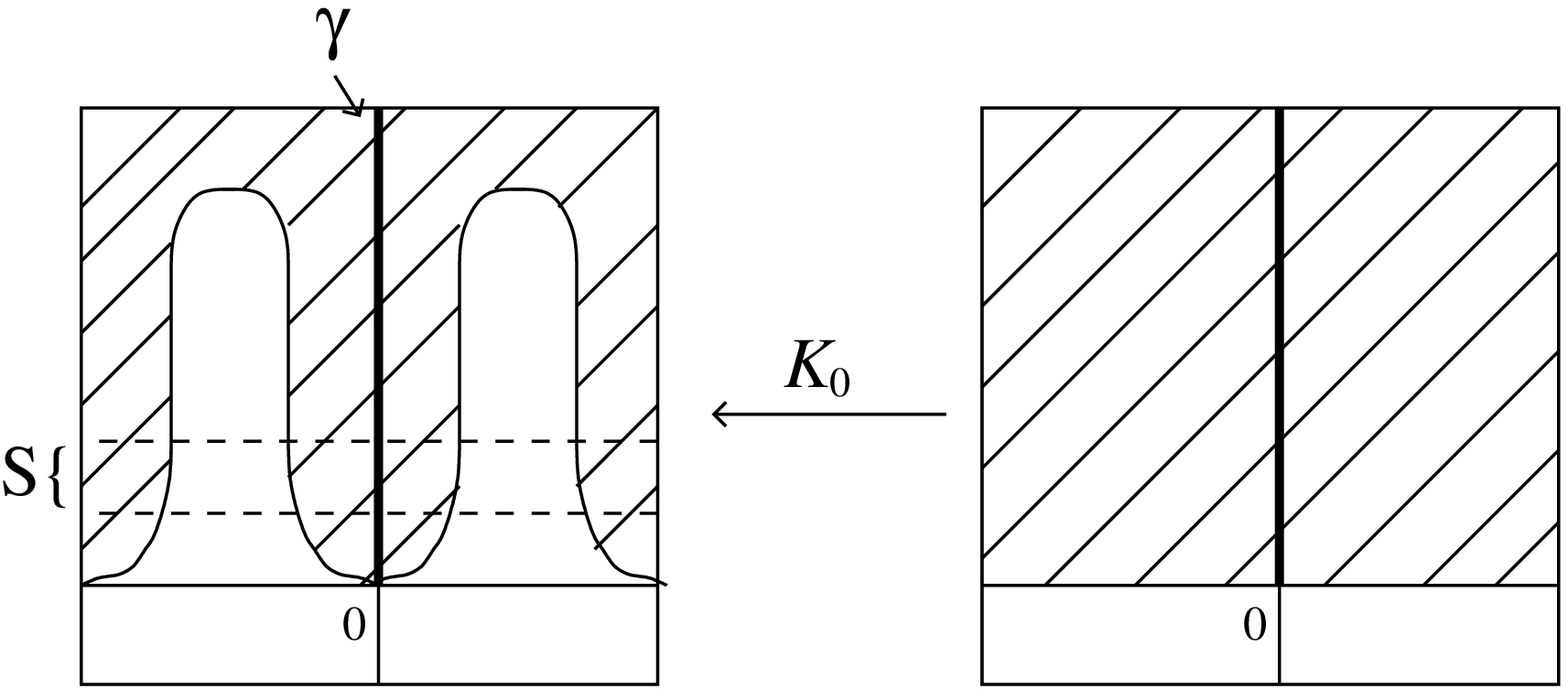,height=4cm}}}}

At stage $n = 0$, $a_0$ may be taken to be any positive integer greater than $5$ and $h_0$ any real greater than
$2$. Only the map $K_0$ has been chosen, and we have only one curve $\gamma := K_0(\xi_0) = \xi_0$, as in the
figure above. This gives the curve $\Gamma := \gamma \cap S$ within the strip $S$.

\medskip

{\hfill {\centerline {\psfig {figure=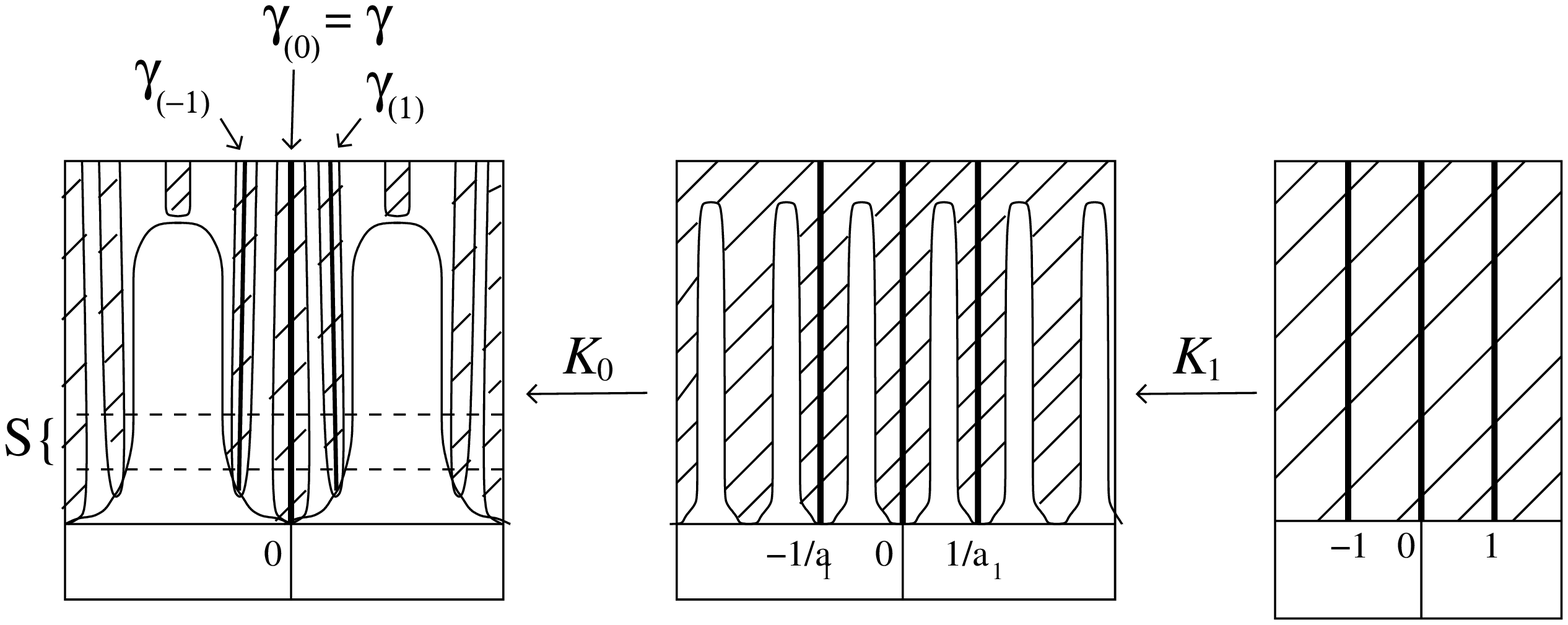,height=4cm}}}}

\medskip

At the next stage $n = 1$, we need to choose the map $K_1$. Now the images of the half-lines $\xi_{-1}, \xi_1$
under $K_1$ are the half-lines $\xi_{-1/a_1}, \xi_{1/a_1}$. So by choosing $a_1$ very large, since $1/a_1$ will be
very small, the curves $K_0(K_1(\xi_{-1})) = K_0(\xi_{-1/a_1}), K_0(K_1(\xi_1)) = K_0(\xi_{1/a_1})$ will be "very
close" to the curve $\gamma = K_0(\xi_0)$. More precisely, these curves will also pass through the strip $S$, be
graphs over the imaginary axis, and if we fix an $r \geq 0$, then the parametrizations with respect to the
imaginary coordinate will be close in $C^r$ norm when $a_1 \to +\infty$. So at stage $n = 1$, we have three
curves, $\gamma_{(-1)} := K_0(\xi_{-1/a_1}), \, \gamma_{(1)} := K_0(\xi_{1/a_1})$ and $\gamma_{(0)} := \gamma$ as
in the figure above, which give the curves $\Gamma_{(-1)} := \gamma_{(-1)} \cap S, \, \Gamma_{(1)} := \gamma_{(1)}
\cap S$ and $\Gamma_{(0)} := \gamma_{(0)} \cap S$ within the strip $S$.

\medskip

{\hfill {\centerline {\psfig {figure=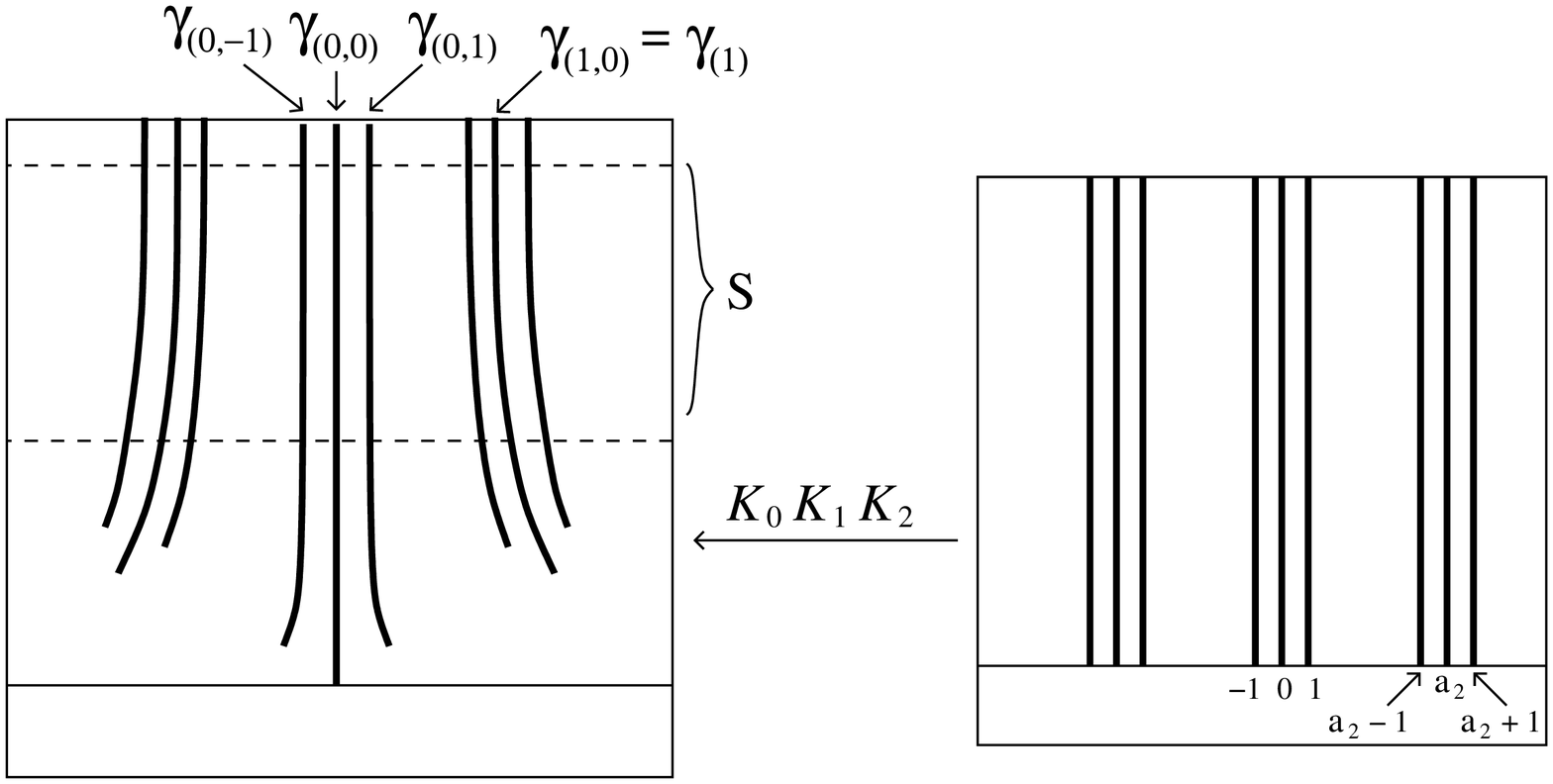,height=5cm}}}}

\medskip

The next stage $n = 2$ is similar. Under the map $K_2$, for $\epsilon_0, \epsilon_1 \in \{ -1, 0, 1 \}$ the
half-lines $\xi_{ \epsilon_0 a_2 + \epsilon_1 }$ are taken to the half-lines $\xi_{ \epsilon_0 + \epsilon_1 / a_2
}$. Thus as in the previous stage, by choosing $a_2$ very large, since $\epsilon_1 / a_2$ will be very small, the
curves $K_0(K_1(K_2(\xi_{ \epsilon_0 a_2 + \epsilon_1 }))) = K_0(K_1(\xi_{ \epsilon_0 + \epsilon_1 / a_2 }))$ will
be "very close" to the curves $K_0(K_1(\xi_{ \epsilon_0 }))$ from the previous stage. As before, all the curves
will be graphs over the imaginary axis which pass through the strip $S$, with each new curve formed in this stage
being close in $C^{r+1}$ norm to a curve from the previous stage. At this stage $n = 2$ we get $9 = 3^2$ curves,
$\gamma_{(\epsilon_0, \epsilon_1)} := K_0 K_1 (\xi_{ \epsilon_0 + \epsilon_1 /a_2})$, such that
$\gamma_{(\epsilon_0, 0)} = \gamma_{(\epsilon_0)}$, for $\epsilon_0, \epsilon_1 \in \{ -1, 0 ,1 \}$. Each gives a
curve $\Gamma_{(\epsilon_0, \epsilon_1)} := \gamma_{(\epsilon_0, \epsilon_1)} \cap S$ within the strip $S$.

\medskip

Proceeding similarly, at each stage $n$, the map $K_n$ is chosen, and we obtain $3^n$ curves indexed by
$n$-tuples of the form $(\epsilon_0, \, \dots, \, \epsilon_{n-1}), \epsilon_j \in \{ -1,0,1 \}, j=0, \dots,
n-1$. Each curve is the image under $K_0 \dots K_{n-1}$ of a half-line $\xi_x$, for $x$ of the form
$$
x(\epsilon_0, \, \dots, \, \epsilon_{n-1}) := \epsilon_0 a_2 \dots a_{n-1} + \epsilon_1 a_3 \dots a_{n-1} +
\dots + \epsilon_{n-2} + {\epsilon_{n-1} \over a_n}
$$
The curves are given by
$$\displaylines{
\gamma_{(\epsilon_0, \dots, \epsilon_{n-1})} := K_0 K_1 \dots K_{n-1} ( \xi_{x(\epsilon_0, \, \dots, \,
\epsilon_{n-1})} ) \quad \hbox{ and satisfy} \cr
\gamma_{(\epsilon_0, \dots, \epsilon_{n-2}, \, 0)} = \gamma_{(\epsilon_0, \dots, \epsilon_{n-2})} \, , \hbox{
where } \epsilon_0, \dots, \epsilon_{n-2}, \epsilon_{n-1} \in \{ -1,0,1 \} \cr
}$$
As before, the curves are analytic, graphs over the imaginary axis and pass through the strip $S$. Thus each
curve
$$
\Gamma_{(\epsilon_0, \dots, \epsilon_{n-1})} := \gamma_{(\epsilon_0, \dots, \epsilon_{n-1})} \cap S
$$
admits a (unique) real-analytic parametrization of the form
$$\eqalign{
\sigma_{(\epsilon_0,\dots,\epsilon_{n-1})} : & \ [ 1,2] \rightarrow
\Gamma_{(\epsilon_0,\dots,\epsilon_{n-1})} \cr
                                             & \ \ \ t \quad \mapsto s_{(\epsilon_0,\dots,\epsilon_{n-1})}(t)
+ it  \cr }$$ This definition is made conditionally since the existence of these parametrizations depends on the
choice of the parameters $(a_n),(h_n)$; that these parameters can indeed be chosen appropriately will be justified
in the following section.

\medskip

As explained above, choosing the $a_n$'s growing fast enough, each parametrization
$\sigma_{(\epsilon_0,\dots,\epsilon_{n-1})}$ will be very close to $\sigma_{(\epsilon_0,\dots,\epsilon_{n-2})}$.
By making them very close in $C^n$ norm say, we can ensure that for any (infinite) sequence $\tilde{\epsilon} =
(\epsilon_0, \epsilon_1, \dots) \in \{ -1,0,1 \}^{{\bf N}}$, the sequence
$(\sigma_{(\epsilon_0,\dots,\epsilon_n)})_{n \geq 0}$ is Cauchy in every $C^r$ norm, and hence converges to a
$C^{\infty}$ parametrization, which we denote by $\sigma_{\tilde{\epsilon}}$, of a smooth curve
$\Gamma_{\tilde{\epsilon}}$. We note that for finite sequences $\tilde{\epsilon} = (\epsilon_0, \dots,
\epsilon_{n-1}, 0,0, \dots)$, the curves $\Gamma_{\tilde{\epsilon}}$ coincide with the curves
$\Gamma_{(\epsilon_0, \dots, \epsilon_{n-1})}$ defined in the finite stages of the construction.

\medskip

We define the set (which will be our comb) $$ {\cal C} := \bigcup_{\tilde{\epsilon} \in \{ -1,0,1 \}^{{\bf N}}}
\Gamma_{\tilde{\epsilon}} $$ (this definition is also made conditionally since the existence of the curves
$\Gamma_{\tilde{\epsilon}}$ for infinite sequences ${\tilde{\epsilon}}$ depends on the choice of the parameters
$(a_n),(h_n)$); again, justification will be provided in the following section). The set ${\cal E}$ which will be
the comb of the main theorem is given by ${\cal E} = E({\cal C})$ (where $E(z) = e^{2\pi i z}$). That these sets
are indeed combs and contained in the hedgehog ${\cal H}$ will be shown in section 3.

\medskip

We note that points of ${\cal C}$ can be represented naturally by pairs $(\tilde{\epsilon},t)$, thus exhibiting
for ${\cal C}$ a product structure of the form Cantor set $\times$ interval, with $\tilde{\epsilon}$ belonging to
the abstract Cantor set $\{ -1,0,1 \}^{{\bf N}}$ with the product topology, and $t$ to the interval $[1,2]$.

\medskip

To study the regularity of the comb ${\cal C}$ however, we will parametrize it by the metric comb $C_0 \times
[1,2]$, where the Cantor set $C_0 \subseteq {\bf R}$ is given by embedding $\{ -1,0,1 \}^{{\bf N}}$ into ${\bf
R}$, via the map $$\eqalign{ \Theta : & \ \{ -1,0,1 \}^{{\bf N}} \, \rightarrow {\bf R} \cr
         & \qquad \tilde{\epsilon} \qquad \quad \mapsto \theta(\tilde{\epsilon}) := {\epsilon_0 \over a_0
a_1} + {\epsilon_1 \over a_0 a_1 a_2} + \dots + {\epsilon_n \over a_0 \dots a_{n+1}} + \dots \cr
}$$
(the map is a homeomorphism onto its image $C_0$ if the $a_n$'s are large enough, for example if $a_n \geq 4$
for all $n$). For sequences $\tilde{\epsilon} = (\epsilon_0, \dots, \epsilon_n,0,0,0,\dots)$ with finitely
many non-zero terms we will denote $\theta(\tilde{\epsilon})$ by $\theta(\epsilon_0,\dots,\epsilon_n)$.

\medskip

The homeomorphism $\Phi : C_0 \times [1,2] \to {\cal C}$ is then given by $$\eqalign{ \Phi : & \ C_0 \times [1,2]
\rightarrow {\cal C} \cr
       & \ (\theta(\tilde{\epsilon}) , \, t) \quad \mapsto  \sigma_{\tilde{\epsilon}}(t) \cr
}$$ (the homeomorphism $\Psi$ in the main theorem is given by composing with the exponential, $\Psi = E \circ
\Phi$). We remark again that this definition is conditional, depending on the choice of parameters $a_n, h_n$; in
the following section 2.3, we show how to choose the parameters so that this map $\Phi$ is well-defined, a
homeomorphism, and satisfies the H\"older and Lipschitz conditions stated in the main theorem, which will be
proved in section 3. We note also that, ensuring that the comb is strictly contained in a period strip for the
exponential, such as $\{ \ -1/2 < $ Re $z < 1/2 \ \}$ (this will be achieved by proper choice of the parameters),
it is not hard to check that these properties for $\Phi$ imply the same for $\Psi$.

\bigskip

\noindent {\bf 2.3  Inductive choice of the parameters.}

\medskip

We let
$$
G_n = u_n + i v_n := K_0 \dots K_n \ , \, n \geq 0
$$

We fix an increasing sequence $\alpha_n \uparrow 1$ ( we will ensure at each stage $n$ in the formation of the
comb that $\Phi$ satisfies an $\alpha_n$-H\"older estimate, so that in the limit $\Phi$ will be $\alpha$-H\"older
for all $\alpha$). Fix a constant $M = 1$ (any $M > 0$ will do); the parameters will be chosen so that the limit
semi-flow ${\cal F}$ is defined on the half-plane ${\bf H}_{-M}$.

\medskip

The following lemma will be used in the construction to ensure that $\Phi$ is injective and $\Phi^{-1}$ Lipschitz:

\medskip

\noindent {\bf Lemma 2.3.1.} \ {\it Given fixed constants $v,h>0$ such that $0<v<h$, for $a \geq 1$ and $x_0 \in
{\bf Z}$, let $z_0 = x_0 + iy_0$, with $y_0$ depending on $a$, be such that Im ${1 \over a} K_{\delta(a,h)}(z_0) =
v$ (note that $y_0$ doesn't depend on $x_0$). Then for $z = x_0+iy$ such that $0 \leq y \leq y_0$, $$ \left|a
\cdot {1 \over a} K_{\delta(a,h)}'(z) \right| = |K_{\delta(a,h)}'(z)| \to \infty \ \hbox{ as } a \to \infty. $$
uniformly in $x_0, \, y$. }

\medskip

\noindent {\bf Proof :} \ We first note that $$\eqalign{ K'(z) & = {e^{2\pi iz} \over -\log(1-e^{2\pi iz})} \cdot
{1 \over (1 - e^{2\pi iz})} \cr
      & = {1 - e^{-e^{2\pi iK(z)}} \over e^{2\pi iK(z)} \cdot e^{-e^{2\pi iK(z)}} } \cr
      & = {e^{e^{2\pi iK(z)}} - 1 \over e^{2\pi iK(z)}} \cr
}$$
Now
$$\eqalign{
{1 \over a}K_\delta(z_0) & = {1 \over a}\left( K(z_0 + i\delta) + iC_\delta \right) \cr
                         & \hbox{ where } C_\delta
\in {\bf R} \hbox{ is a constant depending on } \delta \cr
                         & = {1 \over a}\left(x_0 + K(i(y_0+\delta)) + iC_\delta \right) \cr
}$$
so considering the imaginary part of the above expression gives
$$
{1 \over a}\left( K(i(y_0+\delta)) + iC_\delta \right) = iv \, ,
$$
so
$$\eqalign{
K(i(y_0 + \delta)) & = i(a v - C_\delta) \cr
                   & = i(a v - {1 \over 2\pi}\log\left(-\log(1 - e^{-2\pi \delta})\right)) \cr
                   & = i(a v - {1 \over 2\pi}\log\left( e^{2\pi a h} \log(1+e^{-2\pi \delta}) \right))
\cr
                   & = i(a v - a h - {1 \over 2\pi}\log \log(1 + e^{-2\pi \delta})) \cr
                   & = i(a(v-h) + O(1)) \ \hbox{ as } a \to \infty \cr
}$$ (we have used above the condition on $\delta$ given in section 2.1.4). Thus, letting $w_0 = e^{2\pi i
K(i(y_0+\delta))}$, we have $$ w_0 = e^{2 \pi (a(h-v) + O(1))} \in {\bf R}_+ \ \hbox{ and } \ w_0 \to +\infty \
\hbox{ as } a \to \infty $$ so $$\eqalign{ \left|a \cdot {1 \over a} K_{\delta(a,h)}'(z_0) \right| & =
|K_{\delta(a,h)}'(z_0)| \cr
                                                   & = |K'(z_0+i\delta)| \ \hbox{ (since } K_\delta(z) =
K(z+i\delta) + \hbox{constant}) \cr
                                                   & = |K'(i(y_0+\delta))| \ \hbox{ (since } K(x_0+z) = K(z)
+ x_0 ) \cr
                                                   & = {e^{w_0} - 1 \over w_0} \cr
                                                   & \to +\infty \ \hbox{ as } a \to \infty \cr
}$$
For $z=x_0+iy, \, 0 \leq y \leq y_0$, we note that Im $K(i(y+\delta)) \leq $ Im $K(i(y_0+\delta))$, so $w =
e^{2\pi i K(i(y+\delta))} \geq w_0$, and
$$
\left|a \cdot {1 \over a} K_{\delta(a,h)}'(z_0) \right| = {e^w - 1 \over w} \geq {e^{w_0} - 1 \over w_0} \ ,
$$
a lower bound independent of $x_0, \, y$, hence the convergence is uniform in $x_0, \, y$. \qquad
$\diamondsuit$.

\bigskip

We can now begin the construction. At stage

\medskip

\noindent {\bf n = 0} :  \quad Let $h_0 = 10$. We choose $a_0 \geq 5$ to satisfy the following:

\medskip

Applying Lemma 2.3.1 above to $K_0 = {1 \over a_0} K_{\delta(a_0,h_0)}$, with $v = h_0 -1/2 < h = h_0$ and $a =
a_0$, we require that $a_0$ is chosen large enough so that $$ |a_0 \cdot K_0'(z)| > 1 \ \hbox{ for all } z = x+iy,
\, x \in {\bf Z}, \, 0 \leq y \leq Y_1 $$ where $Y_1$ depending on $a_0$ is chosen such that Im $K_0(iY_1) = h_0 -
1/2$.

\medskip

The parameters $a_0$ and $h_0$ having been chosen, the map $K_0$ is now fixed. We need to choose
additional parameters so that the induction can continue:

\medskip

Using property {\bf (4)} of the maps ${1 \over a} K_{\delta}$, we choose $h_1 > 2$ such that Im $K_0(x+i(h_1-1)) >
h_0 + 1$ for all $x \in {\bf R}$. We choose $h_1' >0$ such that $2 < $ Im $K_0(0+ih_1')$. Also, we have $0 = $ Im
$K_0(0+i0) < 1$; so by continuity, we can choose $\tau_1 > 0$ such that $$ 2 < \hbox{ Im } K_0(x+ih_1') \ , \
\hbox{Im } K_0(x+i0) < 1 \ \hbox{ for } |x| \leq \tau_1 $$ By property {\bf (3)} of the maps ${1 \over a}
K_{\delta}, \ {\partial v_0 \over \partial y} > 0$ on the vertical segment $[0+i0, 0+ih_1']$; thus by continuity
and compactness of the segment, we also take $\tau_1$ small enough to have $$ {\partial v_0 \over \partial y}
(x+iy) > 0 \ \hbox{ for } |x| \leq \tau_1, \, 0 \leq y \leq h_1' $$ The reason for imposing the conditions above
is to guarantee that for $|x| \leq \tau_1$ the curves $K_0(\xi_x)$ pass through the strip $S$ and, in $S$, are
graphs over the imaginary axis, so that in the next stage, we can choose $a_1$ such that $1/a_1 < \tau_1$ to
ensure the same is true for the curves $\gamma_{(-1)} = K_0(\xi_{-1/a_1}), \, \gamma_{(1)} = K_0(\xi_{1/a_1})$.

\medskip

By continuity and compactness of the segment $[0+i0, 0+iY_1]$, we also take $\tau_1$ small enough so that $$ |a_0
\cdot K_0'(z)| > 1 \ \hbox{ for all } z = x+iy, \, |x| \leq \tau_1, \, 0 \leq y \leq Y_1 $$ The set $K_0(\{ \ z =
x+ iy : |x| \leq \tau_1, \, 0 \leq y \leq Y_1 \ \})$ is a neighbourhood of the compact segment $[0+1i, 0+2i]$
(which is the first curve $\Gamma = \gamma \cap S$ occurring in the construction as described in section 2.2). So
we can choose an open rectangle $R$ such that $$ [0+1i, 0+2i] = \Gamma \subset R \subset K_0(\{ \ z = x+ iy : |x|
\leq \tau_1, \, 0 \leq y \leq Y_1 \ \}) $$ We define the domains $D_1$ and $U$ by $$ U := K_0^{-1}(R), \ D_1 = U.
$$ We observe that $\Gamma \subset K_0(U), \ K_0(U) = R$, and $U \subset \{ \ z = x+ iy : |x| \leq \tau_1, \, 0
\leq y \leq Y_1 \ \}$ so that $|K_0'(z)| > 1/a_0$ in $U$.

\bigskip

\noindent ${\bf n \geq 1}$ : \quad We assume the parameters $a_0, \, \dots, \, a_{n-1}, \, h_0, \, \dots, \,
h_{n-1}$ have been chosen, so the maps $K_0, \, \dots , \, K_{n-1}$ are fixed. As induction hypothesis we also
assume that:

\medskip

\noindent {\bf 1.} $h_n > 2$ has been chosen such that Im $K_0 \dots K_{n-1}(x+i(h_n-1)) = v_{n-1}(x+i(h_n-1))>
h_0 + n$ for all $x \in {\bf R}$.

(This condition will allow us to ensure that all elements of the semi-flow ${\cal F}_{0,n}$ have a common periodic
orbit at a height at least $h_0 + n$, so that elements of the limit flow will have periodic orbits above any given
height and hence be nonlinearisable for irrational rotation numbers).

\medskip

\noindent {\bf 2.} Constants $h_n' > 0$ and $\tau_n > 0$ have been chosen such that $$\displaylines{
v_{n-1}(x+i0)<1 \ , v_{n-1}(x+ih_n')>2 \ \hbox{ and}  \cr {\partial v_{n-1} \over \partial y} (x + iy) > 0 \cr
\hbox{ for } |x - x(\epsilon_0,\dots,\epsilon_{n-2},0)| \leq \tau_n, \ 0 \leq y \leq h_n', \cr \hbox{and }
\epsilon_0, \dots, \epsilon_{n-2} \in \{ -1,0,1 \} \ (\hbox{for } n=1, \hbox{ we only have } |x| \leq \tau_1) \cr
}$$ (This condition is to ensure that the images $K_0 \dots K_{n-1}(\xi_x)$ of lines $\xi_x$ lying in the
rectangular strips $\{ \ |x - x(\epsilon_0,\dots,\epsilon_{n-2},0)| \leq \tau_n, \ 0 \leq y \leq h_n' \ \}$ are
'vertical' enough, ie are graphs over the imaginary axis passing through $S$).

\medskip

\noindent {\bf 3.} For all $\epsilon_0, \dots, \epsilon_{n-2} \in \{ -1,0,1 \}$ we have $$ |a_{n-1} \cdot
K_{n-1}'(z)| > 1 \ \hbox{ for all } z = x+iy, \, |x - x(\epsilon_0,\dots,\epsilon_{n-2},0)| \leq \tau_n, \, 0 \leq
y \leq Y_n $$ where $Y_n$ is a constant which has been chosen such that Im $K_{n-1}(Y_n) = h_{n-1} - 1/2$ (for
n=1, the above inequality holds for $|x| \leq \tau_1, \, 0 \leq y \leq Y_1$).

\medskip

\noindent {\bf 4.} Open sets $D_k, \, k=1,\dots,n$ have been chosen, of the form $$ D_k = \bigcup_{
\epsilon_0,\dots,\epsilon_{k-2} \in \{ -1,0,1 \} } U(\epsilon_0,\dots,\epsilon_{k-2}) \ , \ 2 \leq k \leq n $$
(for $k=1$ we only have $D_1 = U$) where each $U(\epsilon_0,\dots,\epsilon_{k-2})$ is a domain such that $$
\Gamma_{(\epsilon_0, \dots, \epsilon_{k-2})} \subset K_0 \dots K_{k-1}(U(\epsilon_0,\dots,\epsilon_{k-2})) $$ (for
$k=1$ we have $\Gamma \subset K_0(U)$), the image of $U(\epsilon_0,\dots,\epsilon_{k-2})$ under $K_{k-1}$ is an
open rectangle $R(\epsilon_0,\dots,\epsilon_{k-2})$, $$ K_{k-1}(U(\epsilon_0,\dots,\epsilon_{k-2})) =
R(\epsilon_0,\dots,\epsilon_{k-2}) $$ (for $k=1$ we have $K_0(U) = R$) and $$ U(\epsilon_0,\dots,\epsilon_{k-2})
\subset \{ \ z = x+iy, \, |x - x(\epsilon_0,\dots,\epsilon_{k-2})| \leq \tau_k, \, 0 \leq y \leq Y_k \ \} $$ (for
$k=1$ we have $U \subset \{ \ z = x+ iy : |x| \leq \tau_1, \, 0 \leq y \leq Y_1 \ \}$) so that $$ |K_{k-1}'(z)| >
1/a_{k-1} \ , \ z \in U(\epsilon_0,\dots,\epsilon_{k-2}) $$ (for $k=1$ we have $|K_0'(z)| > 1/a_0$ in $U$).
Finally, we also assume a nesting condition, $$ K_{k-1}(U(\epsilon_0,\dots,\epsilon_{k-2})) =
R(\epsilon_0,\dots,\epsilon_{k-2}) \subset U(\epsilon_0,\dots,\epsilon_{k-3}) $$ (for $k=1$ this condition is
empty, while for $k=2$ we request $K_1(U(\epsilon_0)) = R(\epsilon_0) \subset U$).

\medskip

These hypotheses imply the following proposition, which will be used
later in proving that $\Phi$ is injective and $\Phi^{-1}$ is Lipschitz:

\medskip

\noindent {\bf Proposition 2.3.2.} \ {\it For $n \geq 1$ and $z, z' \in U(\epsilon_0,\dots,\epsilon_{n-2})$ ($z,
z' \in U$ for $n=1$), we have $$ |K_{n-1}(z) - K_{n-1}(z')| \geq {|z - z'| \over a_{n-1}} $$ }

\medskip

\noindent {\bf Lemma 2.3.3.} \ {\it Let $f:U_1 \to U_2$ be a holomorphic diffeomorphism between two open sets
$U_1, U_2 \subset {\bf C}$ such that $|f'(z)| \geq C > 0$ for all $z \in U_1$, for some constant $C > 0$. Suppose
$U_2$ is convex. Then $$ |f(z) - f(z')| \geq C |z - z'| $$ for all $z,z' \in U_1$. }

\medskip

\noindent {\bf Proof :} Since $U_2$ is convex, given $z, z' \in U_1$ the straight line segment joining $f(z)$ to
$f(z')$ lies in $U_2$, thus applying the Mean value Theorem to $f^{-1}$ gives the desired inequality.
$\diamondsuit$

\medskip

\noindent {\bf Proof of Proposition 2.3.2 :} The result follows from applying the above lemma to the map $K_{n-1}
: U(\epsilon_0,\dots,\epsilon_{n-2}) \to R(\epsilon_0,\dots,\epsilon_{n-2})$, observing that by the induction
hypothesis {\bf 4.}, $$ |K_{n-1}'(z)| > 1/a_{n-1} \ , \ z \in U(\epsilon_0,\dots,\epsilon_{n-2}), $$ and that
$R(\epsilon_0,\dots,\epsilon_{n-2})$, being a rectangle, is convex. $\diamondsuit$

\bigskip

We now outline the conditions required on $a_n$:

\medskip

\noindent {\bf (a)} We require that $a_n \geq 5$ is chosen large enough to satisfy the condition $({\cal C}_n)$,
$$ a_n \geq C_n(a_0, \dots, a_{n-1}, h_0, \dots, h_{n-1}, M) $$ (note $a_n \geq 5$ ensures that $C$ and $C_0$ are
Cantor sets).

\medskip

\noindent {\bf (b)} For the existence of the parametrizations $\sigma_{(\epsilon_0, \dots,\epsilon_{n-1})} : [1,2]
\to \Gamma_{(\epsilon_0,\dots,\epsilon_{n-1})}$ of the curves formed at stage $n$ :

\medskip

We require that $1/a_n < \tau_n$. The induction hypothesis {\bf 2.} on $\tau_n$ and $h_n'$ given
above ensures that then the curves $\gamma_{(\epsilon_0,\dots,\epsilon_n)}$ will pass through $S$, and, in
$S$, be graphs over the imaginary axis. Thus the existence of the parametrizations $\sigma_{(\epsilon_0,
\dots,\epsilon_{n-1})} : [1,2] \to \Gamma_{(\epsilon_0,\dots,\epsilon_{n-1})}$ such that Im
$\sigma_{(\epsilon_0, \dots,\epsilon_{n-1})} (t) = t$ is guaranteed.

\medskip

\noindent {\bf (c)} For the curves formed at stage $n$ to be $C^n$ close to those of stage $n-1$:

\medskip

Since $x(\epsilon_0,\dots,\epsilon_{n-2},\epsilon_{n-1}) = x(\epsilon_0,\dots,\epsilon_{n-2},0) +
\epsilon_{n-1}/a_n \to x(\epsilon_0,\dots,\epsilon_{n-2},0)$ as $a_n \to \infty$, clearly the parametrizations $$
\phi_{(\epsilon_0,\dots,\epsilon_{n-2},\epsilon_{n-1})} : y \mapsto K_0\cdots
K_{n-1}(x(\epsilon_0,\dots,\epsilon_{n-2},\epsilon_{n-1}) + iy) \ , \ 0 \leq y \leq h_n' $$ of the curves
$\gamma_{(\epsilon_0,\dots,\epsilon_{n-1})}$ converge to the parametrizations $$
\phi_{(\epsilon_0,\dots,\epsilon_{n-2},0)} : y \mapsto K_0\cdots K_{n-1}(x(\epsilon_0,\dots,\epsilon_{n-2},0) +
iy) \ , \ 0 \leq y \leq h_n' $$ of the curves $\gamma_{(\epsilon_0,\dots,\epsilon_{n-2},0)}$ in all $C^r$ norms as
$a_n \to \infty$; it follows then from the following standard fact, which we state here without proof, that the
parametrizations $\sigma_{(\epsilon_0,\dots,\epsilon_{n-2},\epsilon_{n-1})}$ also converge in all $C^r$ norms to
the parametrizations $\sigma_{(\epsilon_0,\dots,\epsilon_{n-2},0)}$ as $a_n \to \infty$:

\medskip

\noindent {\bf Lemma 2.3.4.} \ {\it Let $\phi_n : y \mapsto u_n(y) + iv_n(y)$ and $\phi : y \mapsto u(y) + iv(y)$,
defined for $0 \leq y \leq h$, be $C^{\infty}$ parametrizations of smooth curves, such that $v_n'(y) > 0, \, v'(y)
> 0$ for all $n, \, y$, and $v_n(0), v(0) < A < B < v_n(h), v(h)$ for all $n$ for fixed constants $A$ and $B$. We
define the $C^{\infty}$ parametrizations $\sigma_n$, $\sigma$ by $$\displaylines{ \sigma_n(t) = u_n \circ
{v_n}^{-1} (t) + it \ , \ v_n(0) \leq t \leq v_n(h) \cr \sigma(t) = u \circ {v}^{-1} (t) + it \ , \ v(0) \leq t
\leq v(h) \cr }$$ Suppose $||\phi_n - \phi||_{C^r[0,h]} \to 0$ as $n \to \infty$ for all $r \geq 0$. Then
$||\sigma_n - \sigma||_{C^r[A,B]} \to 0$ as $n \to \infty$ for all $r \geq 0$. } It follows that by choosing $a_n$
large enough, we can have $$ ||\sigma_{(\epsilon_0,\dots,\epsilon_{n-2},\epsilon_{n-1})} -
\sigma_{(\epsilon_0,\dots,\epsilon_{n-2})}||_{C^n} \leq {1 \over 2^n} $$ for all
$\epsilon_0,\dots,\epsilon_{n-2},\epsilon_{n-1} \in \{ -1,0,1 \}$.

\medskip

In condition {\bf (e)} below we will require a further control on the $C^0$ distance between these
parametrizations, but for now it follows from the above condition on $a_n$ that for any $\tilde{\epsilon} =
(\epsilon_0,\epsilon_1,\dots) \in {\{ -1,0,1 \}}^{{\bf N}}$, the sequence $\{
\sigma_{(\epsilon_0,\dots,\epsilon_m)} \}_{m\geq 0}$ will be Cauchy in every $C^r$ norm on $[1,2]$, and hence
convergent in every $C^r$ norm to a parametrization $\sigma_{\tilde{\epsilon}}$ of a $C^\infty$ curve
$\Gamma_{\tilde{\epsilon}}$. Thus we are guaranteed that the map $\Phi : (\theta(\tilde{\epsilon}),t) \mapsto
\sigma_{\tilde{\epsilon}}(t)$ is, firstly, well-defined, and secondly, smooth in $t$ for fixed
$\theta(\tilde{\epsilon}) \in C_0$, as stated in the Main Theorem.

\medskip

We state the further conditions required on $a_n$ so that in addition $\Phi$ is a bi-H\"older homeomorphism:

\medskip

\noindent {\bf (d)} The following condition will be used in section 3 to prove that the map $\Phi$ is
$\alpha$-H\"older for every $\alpha, \, 0<\alpha<1$:

\medskip

We would like to choose $a_n$ large enough so that we have an estimate of the form
$$
|\Phi(\theta,t) - \Phi(\theta_0,t)| \leq \left( {1 \over a_0\cdots a_n} \right)^{\alpha_n}
$$
for all $t \in [1,2], \theta_0 = \theta(\epsilon_0,\dots,\epsilon_{n-2}), \theta =
\theta(\epsilon_0,\dots,\epsilon_{n-2},\epsilon_{n-1}) \in C_0$ s.t. $|\theta - \theta_0| = {1 \over
a_0\cdots a_n}$.

\medskip

So fix $t_0 \in [1,2]$, and let $z_0 = x(\epsilon_0,\dots,\epsilon_{n-2},0) + iy_0 = x_0 + iy_0$ and $z =
x(\epsilon_0,\dots,\epsilon_{n-2},\epsilon_{n-1}) + iy = x + iy$ be points s.t. $$\displaylines{ \Phi(\theta_0,
t_0) = K_0\cdots K_{n-1}(z_0), \cr \Phi(\theta,t_0) = K_0\cdots K_{n-1}(z). \cr }$$

\medskip

Now $$ {\partial v_{n-1} \over \partial y} (x_0 + iy_0) > 0 \hbox{ (by the induction hypotheses)}, $$ so by the
Implicit Function Theorem, the equation $$ v_{n-1}(x + iy) = t $$ determines, for $(x,t)$ in a neighbourhood of
$(x_0, t_0)$, $y$ as a smooth function of $(x,t)$ (here we consider $x_0,t_0,y_0$ as fixed, and $x,t$ as variable,
since $x$ depends on $a_n$ which we are varying), say $y = y(x,t)$, such that $y(x_0, t_0) = y_0$. Thus $$
\Phi(\theta,t) = u_{n-1}(x+iy(x,t)) + it $$ is a smooth function of $(x,t)$ in a neighbourhood of $(x_0,t_0)$.
Since this holds for all $(x_0,y_0)$, it follows that $\Phi(\theta,t)$ is a smooth function of $(x,t)$,
$\Phi(\theta,t) = f(x,t)$ say, in a compact neighbourhood of the compact $\{ \ (x,t) : x =
x(\epsilon_0,\dots,\epsilon_{n-2},0), \ \epsilon_0,\dots,\epsilon_{n-2} \in \{ -1,0,1 \}, \ t \in [1,2] \ \}$.
Hence given $t \in [1,2]$ and $\theta_0, \theta \in C_0$ as above with $|\theta - \theta_0| = {1 \over a_0\cdots
a_n}$, we can estimate $$|\Phi(\theta, t) - \Phi(\theta_0, t)| = |f(x,t) - f(x_0,t)| \leq ||f||_{C^1} |x - x_0| =
||f||_{C^1} \cdot {1 \over a_n}$$

We note here that the constant $||f||_{C^1}$ only depends on the maps $K_0, \dots, K_{n-1}$ and the parameters
$a_0, \dots, a_{n-1}$ already chosen, and not on $a_n$. Since $\alpha_n < 1$, it follows that by taking $a_n$
large enough we can ensure that $$\eqalign{ |\Phi(\theta, t) - \Phi(\theta_0, t) & \leq ||f||_{C^1} \cdot {1 \over
a_n } \cr & \leq \left({1  \over a_0\cdots a_{n-1}a_n} \right)^{\alpha_n} \cr }$$

\medskip

\noindent {\bf (e)} The following condition will be used in section 3 to prove that the the map $\Phi$ is
injective and $\Phi^{-1}$ Lipschitz:

\medskip

Applying Lemma 2.3.1 to $K_n = {1 \over a_n} K_{\delta(a_n,h_n)}$, with $v = h_n -1/2 < h = h_n$ and $a = a_n$, we
require that $a_n$ is chosen large enough so that $$ |a_n \cdot K_n'(z)| > 1 \ \hbox{ for all } z = x+iy, \, x \in
{\bf Z}, \, 0 \leq y \leq Y_{n+1} $$ where $Y_{n+1}$ depending on $a_n$ is chosen such that Im $K_n(iY_{n+1}) =
h_n - 1/2$.

\medskip

Since the parametrizations $\sigma_{(\epsilon_0,\dots,\epsilon_{n-2},\epsilon_{n-1})}$ converge uniformly to the
parametrizations $\sigma_{(\epsilon_0,\dots,\epsilon_{n-2},0)}$ when $a_n \to +\infty$, for $a_n$ large enough the
curves $\Gamma_{(\epsilon_0,\dots,\epsilon_{n-2},\epsilon_{n-1})}$ are contained in the neighbourhoods $K_0 \dots
K_{n-1}(U(\epsilon_0,\dots,\epsilon_{n-2}))$ of the curves $\Gamma_{(\epsilon_0,\dots,\epsilon_{n-2})}$. We
require that $a_n$ is chosen large enough so that $$ \Gamma_{(\epsilon_0,\dots,\epsilon_{n-2},\epsilon_{n-1})}
\subset K_0 \dots K_{n-1}(U(\epsilon_0,\dots,\epsilon_{n-2})) $$ for all $\epsilon_0, \dots,
\epsilon_{n-2},\epsilon_{n-1} \in \{ -1,0,1 \}$ (for $n=1$ we require that $\Gamma_{(\epsilon_0)} \subset K_0(U),
\ \epsilon_0 \in \{ -1,0,1 \}$). We also require that
$$||\sigma_{(\epsilon_0,\dots,\epsilon_{n-2},\epsilon_{n-1})} -
\sigma_{(\epsilon_0,\dots,\epsilon_{n-2},0)}||_{C^0} \leq {1/20 \over a_0 \dots a_{n-1} } $$ for all $\epsilon_0,
\dots, \epsilon_{n-1} \in \{ -1,0,1 \}$ (note that $a_0, \dots, a_{n-1}$ have already been chosen and fixed).

\medskip

\noindent {\bf (f) } For the hedgehog to be non-linearisable :

\medskip

Let $P_n$ be the image under $K_n$ of the upper ends of the 1-cylinders in ${\cal S}_{\delta_{n+1}}$, or more
precisely, $$ P_n := z_n + {1 \over a_n} {\bf Z} $$ where $z_n$ is the limit of $K_n(z)$ as Im $z \to +\infty$ in
one of the 1-cylinders of ${\cal S}_{\delta_{n+1}}$ (note $P_n$ is independent of the choice of 1-cylinder, since
their images differ by $1/a_n$ translations), and define the sets $$\eqalign{ O_{n,n+1} & := P_n \cr O_{n-1,n+1} &
:= K_{n-1}(P_n) \cr O_{n-2,n+1} & := K_{n-2} K_{n-1} (P_n) \cr & \dots \cr O_{0,n+1} & := K_0 \dots K_{n-1} (P_n)
\cr }$$ We have:

\medskip

\noindent {\bf Proposition 2.3.5.} \ {\it For all $0 \leq i \leq n+1, m \geq n+1$, for every $F \in {\cal
F}_{i,m}$, the set $O_{i,n+1}$ is a union of periodic orbits (modulo ${\bf Z}$).}

\medskip

\noindent {\bf Proof :} Fix $m \geq n+1$. Consider an  $F \in {\cal F}_{n, m}$. The map $F$ is the conjugate under
$K_n$ of a lift $\tilde{F}_t$ of a map in ${\cal F}_{n+1,m}$ for a time $t  \in {\bf Z} + [-2{a_{n+1}}^{-1},
2{a_{n+1}}^{-1}]$, say $t = l + \alpha, l \in {\bf Z}, \alpha \in [-2{a_{n+1}}^{-1}, 2{a_{n+1}}^{-1}]$. Suppose Im
$z \to +\infty$ in the 1-cylinder of ${\cal S}_{\delta_{n+1}}$ whose ramification point lies at a point $q \in
{\bf Z}$ in the $0$-sheet; then $\tilde{F}_t(z)$ escapes to infinity through the 1-cylinder whose ramification
point lies at the point $q + l$ in the $0$-sheet, and Im $\tilde{F}_t(z) \to +\infty$. Since $P_n$ consists of the
images under $K_n$ of the upper ends of the 1-cylinders, it follows that the orbits under $F$ of points of
$O_{n,n+1} = P_n$ are periodic (modulo ${\bf Z}$) and contained in $O_{n,n+1}$.

\medskip

For $0 \leq i \leq n-1$, since the elements of ${\cal F}_{i,m}$ are given by conjugating those of ${\cal F}_{n,m}$
by the map $K_i \dots K_{n-1}$, and $O_{i,n+1}$ is given by transporting $O_{n,n+1}$ by the same map, the desired
conclusion follows. $\diamondsuit$.

\medskip

Using property {\bf (5) } of the maps ${1 \over a} K_{\delta}$, we choose $a_n$ large enough so that $$ P_n
\subseteq \{ \hbox{ Im } z \geq h_n - 1 \ \}. $$ The induction hypothesis {\bf 1.} then guarantees
 that $O_{0,n+1}$, the image of $P_n$ under $K_0 \dots K_{n-1}$, which is a set of
periodic orbits (modulo ${\bf Z}$) for all the semi-flows ${\cal F}_{0,m}, \, m \geq n+1$, will be at a height at
least $h_0 + n$ in the upper half-plane, so the semi-flows $({\cal F}_{0,k})_{k \geq 0}$ will have non-linearities
at increasing heights.

\bigskip

We now choose and fix $a_n \geq 5$ large enough to satisfy all the conditions {\bf (a)-(f)} above. Thus the map
$K_n$ is determined.

\bigskip

To complete the induction step, it remains to choose the constants $h_{n+1}, \, h_{n+1}', \, \tau_{n+1}$ and
the domain $D_{n+1}$ so that the induction hypotheses are satisfied for the next step.

\medskip

For hypothesis {\bf 1.} :

\medskip

Now $K_n$ is fixed, so using property {\bf (4)} of the maps ${1 \over a}
K_{\delta}$, we choose $h_{n+1} > 2$ such that Im $K_0 \dots K_n(x+i(h_{n+1}-1)) = v_n(x+i(h_{n+1}-1))>
h_0 + (n+1)$.

\medskip

For hypothesis {\bf 2.} :

\medskip

Using property {\bf (4)} of the maps ${1 \over a} K_{\delta}$, we choose $h_{n+1}'
> 0$ such that
$$
\hbox{Im }K_n(x(\epsilon_0,\dots,\epsilon_{n-2},\epsilon_{n-1},0) + ih_{n+1}') = h_n'
$$
(note $x(\epsilon_0,\dots,\epsilon_{n-2},\epsilon_{n-1},0)$'s are integers).

\medskip

It follows from the induction hypothesis {\bf 2.} and condition {\bf (b)} on $a_n$ that $$\displaylines{
v_n(x(\epsilon_0,\dots,\epsilon_{n-2},\epsilon_{n-1},0)+i0) =
v_{n-1}(x(\epsilon_0,\dots,\epsilon_{n-2},\epsilon_{n-1}) + i0) < 1 \ , \cr
v_n(x(\epsilon_0,\dots,\epsilon_{n-2},\epsilon_{n-1},0)+ih_{n+1}') =
v_{n-1}(x(\epsilon_0,\dots,\epsilon_{n-2},\epsilon_{n-1}) + ih_n') > 2 \hbox{ and} \cr {\partial v_n \over
\partial y} > 0 \cr \hbox{ on all segments } [x(\epsilon_0,\dots,\epsilon_{n-2},\epsilon_{n-1},0) + i0,
x(\epsilon_0,\dots,\epsilon_{n-2},\epsilon_{n-1},0) + ih_{n+1}' ] \cr }$$ by the Chain rule, since the
$x(\epsilon_0,\dots,\epsilon_{n-2},\epsilon_{n-1},0)$'s are integers, so on these vertical segments we have
$K_n'(z)$ is real and positive (by properties {\bf (2), (3)} of the maps ${1 \over a} K_{\delta})$, and hence $$
{\partial v_n \over \partial y} = {\partial v_{n-1} \over \partial y} \cdot
K_n'(x(\epsilon_0,\dots,\epsilon_{n-2},\epsilon_{n-1},0) + iy) > 0. $$

\medskip

So by continuity and compactness of the above line segments, we choose $\tau_{n+1} > 0$ small enough such
that
$$\displaylines{
v_n(x+i0)<1 \ , v_n(x+ih_{n+1}')>2 \ \hbox{ and}  \cr
{\partial v_n \over \partial y} (x + iy) > 0 \cr
\hbox{ for } |x - x(\epsilon_0,\dots,\epsilon_{n-2},\epsilon_{n-1},0)| \leq \tau_{n+1}, \ 0 \leq y \leq
h_{n+1}', \cr
\hbox{and } \epsilon_0, \dots, \epsilon_{n-2}, \epsilon_{n-1} \in \{ -1,0,1 \}. \cr
}$$

\medskip

For hypothesis {\bf 3.} :

\medskip

Using condition {\bf (e)} above on $a_n$, compactness of the segments
$[x(\epsilon_0,\dots,\epsilon_{n-2},\epsilon_{n-1},0) + i0,
x(\epsilon_0,\dots,\epsilon_{n-2},\epsilon_{n-1},0) + iY_{n+1}]$ and continuity, we also choose
$\tau_{n+1}$ small enough so that for all $\epsilon_0, \dots, \epsilon_{n-2}, \epsilon_{n-1} \in \{ -1,0,1
\}$ we have
$$
|a_n \cdot K_n'(z)| > 1 \ \hbox{ for all } z = x+iy, \, |x -
x(\epsilon_0,\dots,\epsilon_{n-2},\epsilon_{n-1},0)| \leq \tau_{n+1}, \, 0 \leq y \leq Y_{n+1}
$$

\medskip

For hypothesis {\bf 4.} :

\medskip

We know from condition {\bf (e)} above on $a_n$ that
$$
\Gamma_{(\epsilon_0,\dots,\epsilon_{n-2},\epsilon_{n-1})} \subset K_0 \dots
K_{n-1}(U(\epsilon_0,\dots,\epsilon_{n-2}))
$$
for all $\epsilon_0, \dots, \epsilon_{n-2},\epsilon_{n-1} \in \{ -1,0,1 \}$,
which implies that
$$
(K_0 \dots K_{n-1})^{-1} \left(\Gamma_{(\epsilon_0,\dots,\epsilon_{n-2},\epsilon_{n-1})} \right) \subset
U(\epsilon_0,\dots,\epsilon_{n-2}).
$$
In addition we have
$$
(K_0 \dots K_{n-1})^{-1} \left(\Gamma_{(\epsilon_0,\dots,\epsilon_{n-1})} \right) \subset K_n(
\{ \ z
= x+iy, \, |x - x(\epsilon_0,\dots,\epsilon_{n-1},0)| \leq \tau_{n+1}, \, 0 \leq y \leq
Y_{n+1} \} ).
$$
Since each set $(K_0 \dots K_{n-1})^{-1} \left(\Gamma_{(\epsilon_0,\dots,\epsilon_{n-2},\epsilon_{n-1})}
\right)$ is in fact a compact vertical line segment (contained in the half-line
$\xi_{x_{(\epsilon_0,\dots,\epsilon_{n-2},\epsilon_{n-1})}}$), we can choose open rectangles
$R(\epsilon_0,\dots,\epsilon_{n-2},\epsilon_{n-1})$ such that
$$
(K_0 \dots K_{n-1})^{-1} \left(\Gamma_{(\epsilon_0,\dots,\epsilon_{n-2},\epsilon_{n-1})} \right) \subset
R(\epsilon_0,\dots,\epsilon_{n-2},\epsilon_{n-1})
$$
and
$$\eqalign{
R(\epsilon_0,\dots,\epsilon_{n-2},\epsilon_{n-1}) & \subset U(\epsilon_0,\dots,\epsilon_{n-2}), \cr
R(\epsilon_0,\dots,\epsilon_{n-2},\epsilon_{n-1}) & \subset K_n( \{ \ z= x+iy, \, |x -
x(\epsilon_0,\dots,\epsilon_{n-2},\epsilon_{n-1},0)| \leq \tau_{n+1}, \, 0 \leq y \leq Y_{n+1} \} ). \cr
}$$
We define $U(\epsilon_0,\dots,\epsilon_{n-2},\epsilon_{n-1})$ by
$$
U(\epsilon_0,\dots,\epsilon_{n-2},\epsilon_{n-1}) :=
K_n^{-1}(R(\epsilon_0,\dots,\epsilon_{n-2},\epsilon_{n-1}))
$$
and $D_{n+1}$ by
$$
D_{n+1} := \bigcup_{\epsilon_0,\dots,\epsilon_{n-2},\epsilon_{n-1} \in \{ -1,0,1 \} }
U(\epsilon_0,\dots,\epsilon_{n-2},\epsilon_{n-1}).
$$
We observe that with these definitions all the conditions
of the induction hypothesis {\bf 4.} are satisfied.

\bigskip

Thus the induction hypotheses {\bf 1.} to {\bf 4.} required for stage $n+1$ are satisfied and the induction
can continue.

\medskip

This ends the inductive choice of the parameters.

\bigskip

\noindent {\bf 3.  Proof of the main theorem.}

\medskip

In the course of the induction above we noted that we get a well defined map $\Phi : C_0 \times [1,2] \to {\cal
C}$, which is $C^{\infty}$ smooth in $t \in [1,2]$ for fixed $\theta \in C_0$. It remains to prove that $\Phi$ is
a homeomorphism, $\Phi$ is $\alpha$-H\"older regular for all $\alpha, \, 0 < \alpha < 1$, $\Phi^{-1}$ is
Lipschitz, and the elements $f_t$ of the limit semi-flow are non-linearisable for irrational $t$.

\medskip

We begin with a lemma which will be a basic tool for the proofs which follow.

\medskip

\noindent {\bf Lemma 3.1.} \ {\it Suppose ${\tilde{\epsilon}}^{(1)}, {\tilde{\epsilon}}^{(2)} \in \{ -1,0,1
\}^{\bf N}$ are such that $$ {{\tilde{\epsilon}}^{(1)}}_n =  {{\tilde{\epsilon}}^{(2)}}_n \hbox{ for } n < i_0,
\hbox{ and } {{\tilde{\epsilon}}^{(1)}}_n <  {{\tilde{\epsilon}}^{(2)}}_n \hbox{ for } n = i_0 $$ for some integer
$i_0\geq 0$. Then $$ \theta({\tilde{\epsilon}}^{(1)}) \leq \theta + {1 \over 3} \delta < \theta + {2 \over 3}
\delta \leq \theta({\tilde{\epsilon}}^{(2)}) $$ where $$ \theta =
\theta({\epsilon_0}^{(1)},\dots,{\epsilon_{i_0}}^{(1)},0,0,\dots) \hbox{  and  } \delta = {1 \over a_0 a_1 \dots
a_{i_0 + 1} } . $$ }

\noindent {\bf Proof :}

$$\eqalign{
\theta({\tilde{\epsilon}}^{(1)}) & = { {\epsilon_0}^{(1)} \over a_0 a_1} + \dots + {
{\epsilon_{i_0}}^{(1)} \over a_0\dots a_{i_0+1}} + { {\epsilon_{i_0+1}}^{(1)} \over a_0\dots a_{i_0+1}
a_{i_0+2}} + \dots \cr
        & \leq  \left( { {\epsilon_0}^{(1)} \over a_0 a_1} + \dots + { {\epsilon_{i_0}}^{(1)}
 \over a_0\dots a_{i_0+1}} \right)  + {1 \over a_0\dots a_{i_0+1}} \left( {1 \over a_{i_0+2} } + {1 \over
a_{i_0+2} a_{i_0+3} } + \dots \right) \cr
        & \leq  \theta + \delta \left( {1 \over 4} + {1 \over 4^2} + \dots \right) \cr
        & = \theta  + {1 \over 3} \delta \cr
}$$

Similarly, noting that ${\epsilon_{i_0}}^{(2)} \geq {\epsilon_{i_0}}^{(1)} + 1$, we get
$$
\theta({\tilde{\epsilon}}^{(2)}) \geq  \theta + {2 \over 3} \delta \qquad \diamondsuit
$$

The following lemma will allow us to approximate $\Phi(\theta,t)$,
where $\theta = \theta(\epsilon_0, \epsilon_1, \dots)$, by $\Phi(\theta_n,t)$,
where $\theta_n = \theta(\epsilon_0,\dots,\epsilon_n)$, and hence to
use the estimates obtained at finite stages in the construction in section 2.3.

\medskip

\noindent {\bf Lemma 3.2.} \ {\it Let $\theta = \theta(\epsilon_0, \epsilon_1, \dots)$ and let $\theta_n =
\theta(\epsilon_0,\dots,\epsilon_n)$ for all $n \geq 0$. Then $$ |\Phi(\theta,t) - \Phi(\theta_n,t)| \leq {1/10
\over a_0 \dots a_{n+1}} $$ }

\medskip

\noindent {\bf Proof:} Using the estimate in {\bf (e)} of the inductive construction we have $$\eqalign{
|\Phi(\theta,t) - \Phi(\theta_n,t)| & \leq \sum_{k = n+1}^{\infty} |\Phi(\theta_k,t) - \Phi(\theta_{k-1},t)| \cr &
\leq \sum_{k = n+1}^{\infty} {1/20 \over a_0 \dots a_k } \cr & \leq {1/20 \over a_0 \dots a_{n+1} } \left( 1 + {1
\over 4} + {1 \over 4^2} + \dots \right) \cr & \leq {1/10 \over a_0 \dots a_{n+1}} \cr }$$

We can now prove

\medskip

\noindent {\bf Proposition 3.3.} \ {\it $\Phi$ is $\alpha$-H\"older for every $\alpha, \, 0 < \alpha < 1$.}

We first prove

\noindent {\bf Lemma 3.4.} \ {\it $\Phi$ is $\alpha$-H\"older in $\theta \in C_0$ for every fixed $t \in [1,2]$,
with a uniform H\"older constant independent of $t$.}

\noindent {\bf Proof :} Let $\alpha$ be given, $0 < \alpha < 1$. Choose $N \geq 1$ such that $\alpha_n > \alpha$
for $n \geq N$.

\medskip

Fix $t \in [1,2]$ and let $\theta, \theta' \in C_0, \theta \neq \theta'$, say $\theta = \theta_{\tilde{\epsilon}}$
and $\theta' = \theta_{\tilde{\epsilon'}}$. We note, using condition {\bf (e)} of the construction and the fact
that $\Phi$ takes values in the strip $S$, that $$ |\hbox{Re } \Phi| \leq \sum_{n = 1}^{\infty} {1/20 \over a_0
\dots a_{n-1} } < 1/2, \ 1 \leq \hbox{ Im } \Phi \leq 2 $$ thus $\Phi$ is a bounded function and it suffices to
prove a H\"older estimate on $\Phi$ for all $\theta, \, \theta'$ sufficiently close to each other. So we assume
that ${\tilde{\epsilon}}$ and ${\tilde{\epsilon'}}$ agree up to $(n-1)$ places with $n \geq N$, ie $\epsilon_i =
\epsilon'_i, \, 0 \leq i \leq n-2, \ \epsilon_{n-1} \neq \epsilon'_{n-1}$, say $\epsilon_{n-1} < \epsilon'_{n-1}$.

\medskip

Then $\theta' - \theta \geq {1 \over 3} \delta$, where $\delta = {1 \over a_0\cdots a_n}$.
For $k \geq 0$, we let $\theta_k = \theta(\epsilon_0,\dots,\epsilon_k), \theta'_k =
\theta(\epsilon'_0,\dots,\epsilon'_k)$. Using the estimate in
{\bf (d)} of the inductive construction and Lemma 3.2 above, we have
$$\eqalign{
|\Phi(\theta',t) - \Phi(\theta,t)| & \leq  |\Phi(\theta'_{n-1},t) - \Phi(\theta_{n-1},t)| + |\Phi(\theta',t)
- \Phi(\theta'_{n-1},t)| \cr
         & + |\Phi(\theta,t) - \Phi(\theta_{n-1},t)| \cr
         & \leq |\Phi(\theta'_{n-1},t) - \Phi(\theta_{n-2},t)| + |\Phi(\theta_{n-1},t) - \Phi(\theta_{n-2},t)| \cr
         & + |\Phi(\theta',t) - \Phi(\theta'_{n-1},t)| + |\Phi(\theta,t) - \Phi(\theta_{n-1},t)| \cr
         & \leq {\delta}^{\alpha_n} + {\delta}^{\alpha_n} + \,|\Phi(\theta',t) - \Phi(\theta'_{n-1},t)| + \,|\Phi(\theta,t) -
\Phi(\theta_{n-1},t)|\ \cr
         & \leq 2 {\delta}^{\alpha_n} + {1 \over 10} \delta + {1 \over 10} \delta \cr
         & \leq 3 {\delta}^{\alpha} \cr
         & \leq 3 \cdot 3^{\alpha} \cdot |\theta' - \theta|^{\alpha} \cr
}$$ (note that the H\"older constant appearing here is independent
of $t$) $\diamondsuit$.

\medskip

\noindent {\bf Proof of proposition 3.3 :} We first note that for all ${\tilde{\epsilon}} \in \{ -1,0,1 \}^{{\bf
N}}$, the parametrizations $\sigma_{\tilde{\epsilon}}$ have uniformly bounded derivatives, since $$\eqalign{
||\sigma_{\tilde{\epsilon}}||_{C^1} & \leq ||\sigma_{(0,0,0,\dots)}||_{C^1} + ||\sigma_{\tilde{\epsilon}} -
\sigma_{(0,0,0,\dots)}||_{C^1} \cr
                                    & \leq ||\sigma_{(0,0,0,\dots)}||_{C^1} + \left( {1 \over 2} + {1
\over 2^2} + \dots \right) \cr
                                    & = ||\sigma_{(0,0,0,\dots)}||_{C^1} + 1 = M = \hbox{ constant} \cr
}$$

So, with the same notation as in Lemma 3.4 above, for points $(\theta',t'), \, (\theta,t)$ with $t'$ not
necessarily equal to $t$, using Lemma 3.4 we have
$$\eqalign{
|\Phi(\theta',t') - \Phi(\theta,t)| & \leq |\Phi(\theta',t') - \Phi(\theta,t')| + |\Phi(\theta,t') -
\Phi(\theta,t)| \cr
       & \leq 3 \cdot 3^{\alpha} \cdot |\theta' - \theta|^{\alpha} + ||\sigma_{\tilde{\epsilon}}||_{C^1} |t'
- t| \cr
       & \leq \left( 3 \cdot 3^{\alpha} + M \right) \cdot \left|(\theta',t') - (\theta,t) \right|^{\alpha}
\qquad \diamondsuit \cr
}$$

\noindent {\bf Proposition 3.5.} \ {\it $\Phi$ is injective, and $\Phi^{-1}$ is Lipschitz.}

\medskip

We need the following lemma:

\medskip

\noindent {\bf Lemma 3.6} \ {\it For $z, z' \in U(\epsilon_0,\dots,\epsilon_{n-1},\epsilon_n)$, we have $$ |K_0
\dots K_n(z) - K_0 \dots K_n(z')| \geq {|z - z'| \over a_0 \dots a_n} $$ }

\medskip

\noindent {\bf Proof :} Recalling the nesting condition of the induction hypothesis {\bf 4.}, $$
K_{k+1}(U(\epsilon_0,\dots,\epsilon_{n-1},\epsilon_k)) \subset U(\epsilon_0,\dots,\epsilon_{k-1}) \ , \ k \geq 1
$$ (the condition is $K_1(U(\epsilon_0)) \subset U$ for $k=0$), we can apply Proposition 2.3.2 repeatedly to
obtain $$\eqalign{ |K_0 \dots K_n(z_1) - K_0 \dots K_n(z_2)| & \geq {|K_1 \dots K_n(z_1) - K_1 \dots K_n(z_2)|
\over a_0} \cr
                                          & \geq {|K_2 \dots K_n(z_1) - K_2 \dots K_n(z_2)| \over a_0 a_1}
\cr
                      & \dots \cr
                                          & \geq {|z_1 - z_2| \over a_0 \dots a_n}. \cr
}$$
$\diamondsuit$

\medskip

\noindent {\bf Proof of Proposition 3.5 :} Let $s+it = \Phi(\theta,t), \, s'+it' = \Phi(\theta',t') \in {\cal C}$,
for $(\theta,t), \, (\theta',t') \in C_0 \times [1,2]$, with $\theta = \theta(\tilde{\epsilon})$ and $\theta' =
\theta(\tilde{\epsilon'})$ say, $\theta \neq \theta'$. Let $n \geq 0$ be such that $\epsilon_i = \epsilon'_i$ for
$i \leq n-1$, and $\epsilon_n \neq \epsilon'_n$. Then we have seen before that $|\theta - \theta'| \geq {1 \over
3} \delta$, where $\delta = 1/(a_0 \dots a_{n+1})$; we note also however that $$\eqalign{ |\theta - \theta'| &
\leq {|\epsilon_n - \epsilon'_n| \over a_0 \dots a_{n+1} } + {|\epsilon_{n+1} - \epsilon'_{n+1}| \over a_0 \dots
a_{n+1} a_{n+2} } + \dots \cr
                   & \leq {2 \over a_0 \dots a_{n+1} } \left( 1 + {1 \over 4} + {1 \over 4}^2 + \dots \right)
\cr
                   & \leq {8 \over 3} \delta \cr
}$$
For $k \geq 0$ we let
$$\displaylines{
\theta_k = \theta(\epsilon_0,\dots,\epsilon_{k-1},\epsilon_k), \ \theta'_k =
\theta(\epsilon'_0,\dots,\epsilon'_{k-1},\epsilon'_k) \ , \cr
\Phi(\theta_k,t) = s_k + it, \ \Phi(\theta'_k,t') = s'_k + it'. \cr
}$$
We let $z_n, \, z'_n$ be points such that
$$\eqalign{
K_0 \dots K_n(z_n) & = \Phi(\theta_n,t) = s_n + it \cr
K_0 \dots K_n(z'_n) & = \Phi(\theta'_n,t') = s'_n + it' \cr
}$$
We observe, using condition {\bf (e)} on $a_{n+1}$, that
$$
z_n \in (K_0 \dots K_n)^{-1} \left(\Gamma_{(\epsilon_0,\dots,\epsilon_{n-1},\epsilon_n)} \right) \subset
U(\epsilon_0,\dots,\epsilon_{n-1}).
$$
and
$$
z'_n \in (K_0 \dots K_n)^{-1} \left(\Gamma_{(\epsilon_0,\dots,\epsilon_{n-1},\epsilon'_n)} \right) \subset
U(\epsilon_0,\dots,\epsilon_{n-1}).
$$
Hence applying Lemma 3.6 above gives
$$
|(s_n+it) - (s'_n+it') | = |K_0 \dots K_n(z_n) - K_0 \dots K_n(z'_n)| \geq {|z_n - z'_n| \over a_0 \dots a_n}
$$
Also, by Lemma 3.2 we have
$$\eqalign{
|(s+it) - (s_n+it)| & = |\Phi(\theta,t) - \Phi(\theta_n,t)| \cr
            & \leq {1/10 \over a_0 \dots a_{n+1} } = 1/10 \cdot \delta \cr
}$$
Similarly $|(s'+it') - (s'_n+it')| \leq 1/10 \cdot \delta$. Therefore
$$\eqalign{
|(s+it) - (s'+it')| & \geq |(s_n + it) - (s'_n + it')| - {1 \over 10} \delta - {1 \over 10} \delta \cr
                    & \geq {|z_n - z'_n| \over a_0 \dots a_n} - {1 \over 5} \delta \cr
                    & \geq {|x_n - x'_n| \over a_0 \dots a_n} - {1 \over 5} \delta \cr
                    & \geq {1 \over a_0 \dots a_n a_{n+1}} - {1 \over 5} \delta \cr
                    & = {4 \over 5} \delta \cr
                    & \geq {4 \over 5} \cdot {3 \over 8} |\theta - \theta'| \cr
}$$
Combining this with $|(s+it) - (s'+it')| \geq |t - t'|$ gives
$$\eqalign{
|(\theta,t) - (\theta',t')|  & \leq |\theta - \theta'| + |t -t'| \cr
                             & \leq {5 \over 4} \cdot {8 \over 3} |(s+it) - (s'+it')| + |(s+it) -
(s'+it')| \cr
                             & = {13 \over 3} |(s+it) - (s'+it')| = {13 \over 3} |\Phi(\theta,t) -
\Phi(\theta',t')| \cr }$$ (note that this inequality also holds when $\theta = \theta'$, since then $|(\theta,t) -
(\theta',t')| = |t - t'| \leq {13 \over 3} |t - t'|  = {13 \over 3} |\Phi(\theta,t) - \Phi(\theta',t')|$).

\medskip

Thus $\Phi$ is injective and $\Phi^{-1}$ is Lipschitz. \qquad $\diamondsuit$

\medskip

In particular, $\Phi$ and $\Phi^{-1}$ are continuous, and hence we have

\medskip

\noindent {\bf Corollary 3.7.} \ {\it $\Phi$ is a homeomorphism. Hence the sets ${\cal C}$ and ${\cal E} = E({\cal
C})$ are combs.}

\medskip

\noindent {\bf Remark :} In the course of proving Lemma 3.4, we observed that $|$ Re $\Phi | < 1/2$, so that
${\cal C}$ is contained in the strip $\{ -1/2 < $ Re $z < 1/2 \}$, which we recall was enough to guarantee that
all properties of $\Phi$ also hold for $\Psi$.

\medskip

\noindent {\bf Proposition 3.8.} \ {\it The comb ${\cal C}$ is contained in the hedgehog ${\cal H}$.}

\medskip

\noindent {\bf Proof :} For each curve $\Gamma_{(\epsilon_0, \dots, \epsilon_{n-1})}$, we know that $$
\Gamma_{(\epsilon_0, \dots, \epsilon_{n-1})} \subset K_0 K_1 \dots K_{n-1} ( \xi_{x(\epsilon_0, \, \dots, \,
\epsilon_{n-1})} ) = K_0 K_1 \dots K_{n-1} K_n ( \xi_{a_n \cdot x(\epsilon_0, \, \dots, \, \epsilon_{n-1})} ) $$
Since $a_n \cdot x(\epsilon_0, \, \dots, \, \epsilon_{n-1})$ is an integer, it follows from Proposition 2.2.1 that
each curve $\Gamma_{(\epsilon_0, \dots, \epsilon_{n-1})}$ is contained in ${\cal H}$. Since curves in ${\cal C}$
are uniform limits of such curves and ${\cal H}$ is closed it follows that ${\cal C}$ is contained in ${\cal H}$.
$\diamondsuit$

\medskip

Finally, we observe that the hedgehog is non-linearisable:

\medskip

\noindent {\bf Proposition 3.9.} \ {\it Each $f_t$ of the limit semi-flow is non-linearisable for $t \in C \cap
({\bf R} - {\bf Q})$.}

\medskip

\noindent {\bf Proof :} From Proposition 2.3.5, we know that the set $O_{0,n+1}$ forms a set of periodic orbits
(modulo ${\bf Z}$) for all elements of the semi-flows ${\cal F}_{0,m}$ for all $m \geq n+1$. We observe also that
for all $m \geq n+1$ and all $F \in {\cal F}_{0,m}$, the periods of these orbits are uniformly bounded above by
the cardinal of $O_{0,n+1} / {\bf Z}$ (which is in fact equal to $a_n$). Thus these orbits are also periodic
orbits (modulo ${\bf Z}$) for elements of the limit semi-flow ${\cal F}$, and each $O_{0,n+1}$ forms a set of
periodic orbits (modulo ${\bf Z}$) for all elements of ${\cal F}$. By {\bf (f)} of the inductive construction,
$P_n$ is contained in the half-plane $\{ $ Im $z \geq h_n - 1 \ \}$, and hence by the induction hypothesis {\bf
1.}, $O_{0,n+1}$ is contained in the half-plane $\{ $ Im $z \geq h_0 + n \ \}$. Thus each $f_t$ has a sequence of
periodic orbits $E(O_n)$ which accumulate $0$, and is hence non-linearisable for irrational $t \in C \cap ({\bf R}
- {\bf Q})$. $\diamondsuit$

\medskip

This ends the proof of the main theorem.

\bigskip

\noindent {\bf 4. Remarks.}

\medskip

\noindent 1. In the construction given here, we make the heights of the periodic orbits $O_{0,n+1}$ increase to
infinity, so that the hedgehog is non-linearisable. If instead the heights $h_n$ are chosen differently so that
the heights of the orbits $O_{0,n+1}$ do not escape to infinity but converge to a finite height, then the limit
dynamics is linearisable above this height, and we obtain a linearisable hedgehog containing a Siegel disk (in
fact the construction in [7] is carried out this way). It is possible in this case too to construct smooth combs
inside the hedgehog, and to make sure that they are non-trivial, lying outside the Siegel disk.

\medskip

\noindent 2. We note that though the hedgehog is invariant under the dynamics, the comb itself is not.

\medskip

\noindent 3. Since the comb ${\cal C}$ is bi-H\"older-$\alpha$ equivalent to the comb $C_0 \times [1,2]$ for all
$\alpha, \, 0 < \alpha < 1$, ${\cal C}$ has the same Hausdorff dimension as $C_0 \times [1,2]$, which is the
dimension of $C_0$ plus one. It is easy to see from the explicit representation of $C_0$ in terms of the $a_n$'s
that $C_0$ can be made zero-dimensional by taking the $a_n$'s growing fast enough. Thus we can ensure that the
comb ${\cal C}$ has Hausdorff dimension one. We can also get combs of Hausdorff dimension one by taking small
enough Cantor subsets $C_0' \subset C_0$ of $C_0$, with $C_0'$ zero-dimensional.

\medskip

\noindent 4. The rotation numbers $t \in C$ for the limit dynamics include numbers of the form $$ t = \sum_{n =
0}^{\infty} {\epsilon_n \over a_0 a_1 \dots a_n} $$ where $\epsilon_n \in \{ \, -1,0,1 \, \} \, , \, n \geq 0$. As
the construction requires the $a_n$'s to be growing rapidly, this means that we obtain very Liouville rotation
numbers.

\bigskip

\noindent {\bf 5. Acknowledgements.}

\medskip

The main result of this article, which is based on part of my Ph.D. thesis, was first stated by my advisor
R.Perez-Marco in [7]. I am very grateful to him for proposing this problem to me and also for explaining how to
use his techniques of tube-log Riemann surfaces to solve it. I thank also the referee for many helpful suggestions
which have improved this article considerably.

\bigskip

{\centerline {\bf References}}

\bigskip

\medskip

\noindent [1] G.D.BIRKHOFF, {\it Surface transformations and their dynamical applications}, Acta Mathematica, {\bf
43}, (1920). Also, Collected Mathematical Papers II, p.111

\medskip

\noindent [2] H.CREMER, {\it Zum Zentrumproblem}, Math Ann., vol 98, 1928, p. 151-153

\medskip

\noindent [3] H.CREMER, {\it \"Uber die H\"aufigkeit der Nichtzentren}, Math Ann., vol 115, 1938, p.573-580

\medskip

\noindent [4] R.PEREZ-MARCO, {\it Fixed points and circle maps}, Acta Mathematica, {\bf 179:2}, 1997, p.243-294.

\medskip

\noindent [5] R.PEREZ-MARCO, {\it Topology of Julia sets and hedgehogs}, Preprint, Universit\'e de PARIS-SUD,
1994.

\medskip

\noindent [6] R.PEREZ-MARCO, {\it Uncountable number of symmetries for non-linearisable holomorphic dynamics},
Inventiones Mathematicae, {\bf 119}, {\bf 1}, p.67-127,1995; (C.R. Acad. Sci. Paris, {\bf 313}, 1991, p. 461-464).

\medskip

\noindent [7] R.PEREZ-MARCO, {\it Siegel disks with smooth boundary}, Preprint, Universit\'e de PARIS-SUD, 1997,
www.math.ucla.edu/$\tilde{\ }$ricardo, submitted to the Annals of Mathematics in 1997.

\medskip

\noindent [8] R.PEREZ-MARCO, {\it Hedgehog's Dynamics}, Preprint, UCLA, 1996.

\medskip

\noindent [9] R.PEREZ-MARCO,  {\it Sur les dynamiques holomorphes non lin\'earisables et une conjecture de V.I.
Arnold}, Ann. Scient. Ec. Norm. Sup. 4 serie, {\bf 26}, 1993, p.565-644; (C.R. Acad. Sci. Paris, {\bf 312}, 1991,
p.105-121).

\medskip

\noindent [10] C.L.SIEGEL, {\it Iteration of Analytic Functions}, Ann. Math., vol.43, 1942, p. 807-812

\medskip

\noindent [11] J.C.YOCCOZ, {\it Petits diviseurs en dimension 1}, S.M.F., Asterisque {\bf 231} (1995)

\medskip

Email address : kbiswas@math.ucla.edu

\end